\documentclass{amsart}

\usepackage{amsmath,amsthm,amssymb,amscd,mathtools}
\usepackage{amsfonts,yfonts}
\usepackage{rotating}
\usepackage{euscript}
\usepackage{pst-node}
\usepackage{epsfig,verbatim}
\usepackage{pictexwd,dcpic}
\usepackage{enumerate}
\usepackage{caption}
\usepackage{verbatim}

\def\be{\begin{equation}}
\def\ee{\end{equation}}

\def\C{{\mathbb C}} 
\def\f{\EuScript}
\def\N{{\mathbb N}} 
\def\P{{\mathbb P}}
\def\Z{{\mathbb Z}}
\def\R{{\mathbb R}} 
\def\Q{{\mathbb Q}}

\def\e{\eqref}
\def\phi{{\varphi}}
\def\v{{\varepsilon}} 
\def\tt{\widetilde}
\def\deg{{\rm deg\,}}

\def\GCD{{\rm GCD }}

\def\Aut{{\rm Aut }}

\def\bp{\begin{proposition}}
\def\ep{\end{proposition}}

\def\bt{\begin{theorem}}
\def\et{\end{theorem}}
\def\br{\begin{remark}}
\def\er{\end{remark}}
\def\be{\begin{equation}}
\def\bee{\begin{equation*}}
\def\l{\label}

\def\e{\eqref}

\def\ee{\end{equation}}
\def\eee{\end{equation*}}
\def\bl{\begin{lemma}}
\def\el{\end{lemma}}
\def\bc{\begin{corollary}}
\def\ec{\end{corollary}}
\def\pr{\noindent{\it Proof. }}

\def\bd{\begin{definition}}
\def\ed{\end{definition}}
\def\t{\widetilde}

\def\hat{\widehat}

\def\c{\mathcal}

\newtheorem{theorem}{Theorem}[section]
\newtheorem{lemma}[theorem]{Lemma}
\newtheorem{corollary}[theorem]{Corollary}
\newtheorem{proposition}[theorem]{Proposition}
\theoremstyle{definition}
\newtheorem{remark}[theorem]{Remark}

\mathtoolsset{showonlyrefs}

\begin{document} 
\title{ Invariant curves for endomorphisms of $\P^1\times \P^1$}
\author{Fedor Pakovich}
\thanks{
This research was supported by ISF Grant No. 1432/18}
\address{Department of Mathematics \\
Ben-Gurion University of the Negev \\
P.O.B. 653 Beer Sheva \\
8410501 Israel
}
\email{pakovich@math.bgu.ac.il}

\begin{abstract} Let  $A_1, A_2\in \C(z)$ be rational functions
of degree at least two that are neither Latt\`es maps nor conjugate to $z^{\pm n}$ or $\pm T_n.$
We  describe invariant, periodic, and preperiodic algebraic curves for endomorphisms 
of $(\P^1(\C))^2$ of the form $(z_1,z_2)\rightarrow (A_1(z_1),A_2(z_2)).$ 
In particular, we show that if $A\in \C(z)$ is not a ``generalized Latt\`es map'', then 
any $(A,A)$-invariant curve has genus zero and can be parametrized by rational functions commuting with $A$. As an application, for $A$ defined over a  subfield $K$ of $ \C$  we give a criterion for a point of $(\P^1(K))^2$ to have a Zariski dense $(A, A)$-orbit in terms of canonical heights, and  deduce from this criterion a version of a conjecture of Zhang on the existence of rational points with Zariski dense forward orbits. We also prove a result about functional decompositions of iterates of rational functions, which  implies in particular that there exist at most finitely many $(A_1, A_2)$-invariant curves of any given bi-degree $(d_1,d_2).$

\end{abstract} 

\maketitle

\begin{section}{Introduction}
Let $A$ be a rational function of one complex variable. We say that $A$ is {\it special} if it is either a Latt\`es map, or it is conjugate to $z^{\pm n}$ or $\pm T_n.$  In this paper, we describe invariant and, more generally, periodic and preperiodic algebraic curves for endomorphisms $(A_1,A_2)
:\, (\P^1(\C))^2\rightarrow  (\P^1(\C))^2$ given by the formula
 \be \l{ss} (z_1,z_2)\rightarrow (A_1(z_1),A_2(z_2)),\ee where  $A_1$ and $A_2$ are non-special rational functions of degree at least two. 
Note that describing invariant varieties for more general endomorphisms 
\be \l{s2} (z_1,z_2, \dots z_n)\rightarrow (A_1(z_1),A_2(z_2),  \dots A_n(z_n)), \ \ \ n\geq 2, \ee reduces to describing invariant curves for endomorphisms
\eqref{ss} (see \cite{m}, \cite{ms} and  \cite{gny1}, \cite{k}). On the other hand, an arbitrary dominant endomorphism of $(\P^1(\C))^n$ has the form 
$$ (z_1,z_2, \dots z_n)\rightarrow (A_1(z_{\sigma(1)}),A_2(z_{\sigma(2)}),  \dots A_n(z_{\sigma(n)})) $$
for some permutation $\sigma\in S_n$,
implying that  some of its iterates has form \eqref{s2}.

Invariant curves for endomorphisms \eqref{ss} with {\it polynomial} $A_1$, $A_2$  were studied in the  paper  of Medvedev and Scanlon  \cite{ms}. In particular, it was shown in \cite{ms} that if   $A_1$ and $A_2$  are not 
conjugate to powers $z^{n}$ or Chebyshev polynomials $\pm T_n$, then any irreducible algebraic 
$(A_1,A_2)$-invariant curve has genus zero and can be parametrized by polynomials $X_1,$ $X_2$ satisfying the system of functional equations 
\be \l{sys} A_1\circ X_1=X_1\circ B, \ \ \ \  \ \   A_2\circ X_2=X_2\circ B \ee
for some polynomial $B.$ Using the theory of 
functional decompositions of polynomials developed by Ritt (\cite{rpc}), Medvedev and Scanlon 
investigated system \eqref{sys} in detail and obtained a description of  $(A_1,A_2)$-invariant curves. Specifically, for $A_1=A_2$ the main result of \cite{ms} about invariant curves can be formulated as follows:
if a polynomial $A$  is not conjugate to  $z^{n}$ or $\pm T_n$, then 
any irreducible $(A,A)$-invariant curve is a graph $z_2=X(z_1)$ or $z_1=X(z_2)$, where $X$ is a 
polynomial commuting with $A$. 
The classification of invariant curves obtained by  Medvedev and Scanlon 
has numerous applications in arithmetic dynamics 
(see e. g. \cite{bm}, \cite{fg},   \cite{gn}, \cite{gn1},  \cite{gn2}, \cite{gny2}, \cite{gx}, \cite{n}), 
and the  goal of this paper is to obtain a generalization of this classification to arbitrary non-special rational functions $A_1$ and $A_2$.
%For  non-special rational functions $A_1$ and $A_2$, 
For such functions, any $(A_1,A_2)$-invariant curve still has genus zero and 
  can be parametrized by {\it rational} functions $X_1,$ $X_2$ satisfying 
\eqref{sys} for some {\it rational} function $B.$ In particular, the existence of invariant curves implies the equality $\deg A_1=\deg A_2.$ However, the Ritt theory of polynomial decompositions used in \cite{ms} for the analysis of \eqref{sys} does not extend to rational functions. Furthermore, 
one of the key ingredients of the method of \cite{ms}, the so-called ``first Ritt theorem'', is known not to be true in the rational case (see e. g. \cite{mp}). 
Note that results of \cite{ms} about invariant curves  can be 
proved by a different method, which does not rely on the first Ritt theorem (see \cite{pj}). Nevertheless, the method of \cite{pj} is also restricted to the polynomial case. 

Since rational functions parametrizing invariant curves  for endomorphisms  \eqref{ss} satisfy system \eqref{sys},  the problem of describing invariant curves 
 is closely related to the problem of describing  
 {\it semiconjugate rational functions}, that is, 
 rational solutions of the functional equation 
\be \l{ssee} A\circ X=X\circ B.\ee 
A comprehensive description of solutions of \eqref{ssee}  was obtained  in the series of papers \cite{semi}, \cite{rec},  \cite{dyna}, \cite{lattes},  \cite{fin}, and  
in this paper we apply the main results of  \cite{semi} and \cite{lattes} 
to system \eqref{sys}. 
To formulate our results explicitly we recall several definitions. For the rest of this paper, we use the standing convention that ``rational function'' means ``nonconstant rational function''. 

An {\it orbifold} $\f O$ on $\P^1(\C)$  is a ramification function $\nu:\P^1(\C)\rightarrow \mathbb N$ which takes the value $\nu(z)=1$ except at a finite set of points.  
If $f$ is a rational function and $\f O_1$, $\f O_2$ are orbifolds with ramification functions $\nu_1$ and $\nu_2$, then we say that  
$f:\,  \f O_1\rightarrow \f O_2$ is  {\it a covering map} 
between orbifolds
if for any $z\in \P^1(\C)$ the equality 
$$ \nu_{2}(f(z))=\nu_{1}(z)\deg_zf$$ holds.
In case   
the weaker condition 
$$ \nu_{2}(f(z))=\nu_{1}(z)\GCD(\deg_zf, \nu_{2}(f(z))$$ is satisfied,  
we say that $f:\,  \f O_1\rightarrow \f O_2$ is  {\it a  minimal holomorphic  map} 
between orbifolds.
In these terms, {\it a  Latt\`es map} can be defined as a rational function $A$ of degree at least two such that $A:\f O\rightarrow \f O$ is a  covering self-map
for some orbifold $\f O$ (see \cite{mil2}).  
Following \cite{lattes}, we say 
 that  $A$ is {\it a generalized Latt\`es map}  if there exists an orbifold $\f O$ 
 distinct from the non-ramified sphere such that  $A:\,  \f O\rightarrow \f O$ is  a  minimal holomorphic  map. Note that similar to ordinary  Latt\`es maps, generalized Latt\`es maps can be characterized in terms of semiconjugacies and group actions  (see \cite{lattes}).

Let $A_1$, $A_2$  $X_1,$ $X_2,$ $B$ be 
rational functions such that the diagram
\be \l{ddii}
\begin{CD} 
(\P^1(\C))^2 @>(B,B)>>(\P^1(\C))^2 \\ 
@V (X_1,X_2)  VV @VV  (X_1,X_2) V\\ 
 (\P^1(\C))^2 @>(A_1,A_2)>> (\P^1(\C))^2
\end{CD}
\ee
commutes. Then the image of $\P^1(\C)$ in $(\P^1(\C))^2$ under the map
\be\l{parr}  t\rightarrow (X_1(t),X_2(t))\ee is an $(A_1,A_2)$-invariant algebraic 
curve $\c C $, since the diagonal $\Delta$ in $(\P^1(\C))^2$ is  $(B,B)$-invariant and 
$\c C=(X_1,X_2)(\Delta)$.
For brevity, we say that the map  \eqref{parr} is a {\it paramet\-ri\-zation} of the curve $\c C.$
We emphasize however that  
such a parametrization is not necessarily generically one-to-one, that is, we do not assume that
 $X_1$ and $X_2$  satisfy the condition $\C(X_1,X_2)=\C(z)$.

In like manner, if $A_1$, $A_2$  $Y_1,$ $Y_2,$ $B$ are  
rational functions such that the diagram
\be \l{iidd}
\begin{CD} 
(\P^1(\C))^2 @>(A_1,A_2)>>(\P^1(\C))^2 \\ 
@V (Y_1,Y_2)  VV @VV  (Y_1,Y_2) V\\ 
 (\P^1(\C))^2 @>(B,B)>> (\P^1(\C))^2
\end{CD}
\ee
commutes, then the algebraic
curve $\c E=(Y_1,Y_2)^{-1}(\Delta)$, defined by the equation  $ Y_1(x)-Y_2(y)=0,$ satisfies $(A_1,A_2)(\c E)\subseteq \c E$. Therefore, each component of $\c E$ is $(A_1,A_2)$-preperiodic and at least one of these components is $(A_1,A_2)$-periodic.

Our first result  provides a description of $(A_1,A_2)$-invariant curves in case that  $A_1$ and $A_2$ are not generalized Latt\`es maps through a system of functional equations involving functional decompositions of iterates of 
$A_1$, $A_2$ and 
 diagrams \eqref{ddii}, \eqref{iidd}.

\bt \l{1} Let $A_1$, $A_2$ be rational functions of degree at least two that are not generalized Latt\`es maps, and 
$\c C$ an irreducible algebraic curve  in $(\P^1(\C))^2$ that is not a vertical or horizontal line. Then $\c C$ is $(A_1,A_2)$-invariant  if and only if 
 there exist  rational functions $X_1,$ $X_2,$ $Y_1,$ $Y_2,$ $B$  such that: 
\begin{enumerate} 
\item[1.]
 The diagram 
\be  \l{xx}
\begin{CD} 
(\P^1(\C))^2 @>(B,B)>>(\P^1(\C))^2 \\ 
@V (X_1,X_2)  VV @VV  (X_1,X_2) V\\ 
 (\P^1(\C))^2 @>(A_1,A_2)>> (\P^1(\C))^2
\\ 
@V (Y_1,Y_2)  VV @VV  (Y_1,Y_2) V\\ 
 (\P^1(\C))^2 @>(B,B)>> (\P^1(\C))^2 
\end{CD}
\ee
commutes, \item[2.]  The 
equalities 
\be \l{en1} X_1\circ Y_1=A_1^{\circ d}, \ \ \ \ \ \ \ 
 X_2\circ Y_2=A_2^{\circ d},\ee
 \be\l{en2} Y_1\circ X_1=Y_2\circ X_2=B^{\circ d}\ee hold for some $d\geq 0$, \item[3.]  The map $t\rightarrow (X_1(t),X_2(t))$ is a parametrization of  $\c C$. 
\end{enumerate}
\et
Note  that the top square of \eqref{xx} is obtained from elementary considerations about parametrizations of invariant curves in the same way as in the paper \cite{ms} in the polynomial case. On the other hand, the bottom square is based on results  \cite{semi} and \cite{lattes}, and requires the assumption that $A_1$ and $A_2$ are not generalized Latt\`es maps. 

Let us mention that, among other things, Theorem \ref{1}  implies 
that $\c C$ is  
a component of the ``separate variable'' curve \be \l{giv} \c E: Y_1(x)-Y_2(y)=0.\ee Thus, Theorem \ref{1} provides us both with the parametrization  of  $\c C$ and with the equation of a curve having $\c C$ as a
 component. Moreover, both these characterizations of invariant curves are obtained from decompositions of iterates \eqref{en1} subject to special restrictions.
Note also that condition \eqref{en1} yields that 
$$(A_1,A_2)^{\circ d}(\c E)=\c C,$$ that is, all components of curve \eqref{giv} are eventually mapped to the curve $\c C$.

Theorem \ref{1} permits us to describe 
also  $(A_1,A_2)$-periodic and preperiodic curves.   
Specifically, we show that under the assumptions of Theorem \ref{1} a curve   
$\mathcal C$ is $(A_1,A_2)$-periodic  if and only if 
 there exist  rational functions  $X_1,$ $X_2,$ $Y_1,$ $Y_2$ such that
the equalities  
$$
X_1\circ Y_1=A_1^{\circ d}, \ \ \ \ 
 X_2\circ Y_2=A_2^{\circ d}, \ \ \ \  Y_1\circ X_1=Y_2\circ X_2$$ hold for some $d\geq 0$, and the map  $t\rightarrow (X_1(t),X_2(t))$ is a parametrization of  $\c C$.  On the other hand, a curve   $\c C$ is  $(A_1,A_2)$-preperiodic  if and only if 
 there exist  rational functions as above such that 
$\c C$ is a component of curve \eqref{giv}
(Theo\-rem \ref{2}). 
Finally, we show that describing $(A_1,A_2)$-periodic and preperiodic curves for arbitrary non-special rational functions  $A_1$ and $A_2$ reduces to the case where $A_1$ and $A_2$ are not generalized Latt\`es maps (Theorem \ref{3}).  

In a sense, describing $(A_1,A_2)$-periodic and preperiodic curves  reduces to the case  $A_1=A_2=A$ (see  Corollary \ref{cor1}). For this case,  we give the following alternative description of invariant curves, providing an analogue of the result of Medvedev and Scanlon cited above.

\bt \l{1+} Let $A$ be a rational function of degree at least two that is not a generalized Latt\`es map,  and 
$\c C$ an irreducible algebraic curve  in $(\P^1(\C))^2$  that is not a vertical or horizontal line. Then $\c C$ is $(A,A)$-invariant  if and only if 
 there exist  rational functions  $U_1,$ $U_2,$ $V_1,$ $V_2$    commuting with  $A$ such that 
the equalities 
\be \l{en11} U_1\circ V_1=U_2\circ V_2=A^{\circ d},\ee
 \be\l{en21} V_1\circ U_1=V_2\circ U_2=A^{\circ d}\ee hold for some $d\geq 0$ and   
the map $t\rightarrow (U_1(t),U_2(t))$ is a parametrization of  $\c C$.
\et

As an application of Theorem \ref{1+},  for $A$ defined over a number field $K$   we give a criterion for a point of $(x_0,y_0)\in (\P^1(K))^2$ to have a Zariski dense $(A, A)$-orbit in terms of canonical heights of $x_0$ and $y_0$. 
Let us denote by $h$ the Weil height on $\P^1(K)$ and by $\hat h_A$ the corresponding  canonical height associated to $A$. The simplest examples of points  with non-dense $(A,A)$-orbits are points $(x_0,y_0)$  such that $x_0$ or $y_0$ is 
$A$-preperi\-odic. Further examples are points of the form $(x_0,A^{\circ l}(x_0))$ or $(A^{\circ l}(x_0),x_0)$, where $x_0\in \P^1(K)$ and $l\geq 0,$ since such points belong to the 
curves 
\be \l{kria} A^{\circ l}(x)-y=0, \ \ \ \ x-A^{\circ l}(y)=0,\ee 
which  are  $(A,A)$-invariant.  
The canonical heights of the last kind of points obviously satisfy the relation  
%$$\hat h_A(y_0)=n^l\hat h_A(x_0),\ \ \ \ \ \ \ l\in \Z,$$  
$\hat h_A(y_0)=n^l\hat h_A(x_0),$ where $n=\deg A$ and $l\in \Z,$ 
and our main result about orbits states  that a similar relation is sa\-tisfied for any point  $(x_0,y_0)$  whose $(A,A)$-orbit is not dense, provided that $x_0$ and $y_0$  are not $A$-preperiodic.

\bt \l{dh} Let $K$ be a number filed and $A$ a non-special rational function of degree
 $n\geq 2$ defined over $K$. 
Then the $(A,A)$-orbit of a point $(x_0,y_0)\in (\P^1(\overline{K}))^2$ is Zariski dense in $(\P^1(\C))^2$, unless 
either $x_0$ or $y_0$ is a preperiodic point of $A$, or the canonical heights of $x_0$ and $y_0$ satisfy the condition 
\be \l{osn} \hat h_A(y_0)=n_0^{l}\hat h_A(x_0),\ \ \ \ \ \ \ l\in \Z,\ee where $n_0$ is a minimum natural number such that $n=n_0^k$ for some $k\geq 1.$ 
\et

Using instead of the Weil height the Moriwaki height, we also provide an analogue of Theorem \ref{dh} for an arbitrary subfield $K$ of $\C$ finitely generated over $\Q$  (Theo\-rem \ref{dh++}). This allows us to prove a variant of a conjecture of Zhang (\cite{z}) 
on the existence of Zariski dense orbits for endomorphisms of varieties. Namely, we show that if $K$ is a subfield of $\C$ and $A_1,A_2\in K(z)$ are  non-special rational functions of degree at least two, then
there is a point in $(\P^1( K))^2$ whose $(A_1,A_2)$-forward orbit is Zariski dense in $(\P^1( K))^2$ 
(Theorem \ref{z2}).  For algebraically closed fields, this result was established previously in the appendix to the paper \cite{k} as a corollary of the main result of \cite{k} about the existence of Zariski dense orbits for endomorphisms of projective surfaces. The benefits of our approach are that it does not require $K$ to  be algebraically closed, and it permits construct points  with dense orbits in an effective way. 
%in terms of heights.

Since for any rational function $A\in \C(z)$ and integer $l\geq 0$ the  curves \eqref{kria} are 
$(A,A)$-invariant, one cannot expect to bound the total number of $(A_1,A_2)$-invariant curves.
Nevertheless, we show (Theorem \ref{4}) that for any  rational functions $A_1$, $A_2$ of degree $m\geq 2$ there exist at most finitely many $(A_1, A_2)$-invariant curves of any given bi-degree $(d_1,d_2)$, and that the number of such curves can be bounded in terms of $d_1,$ $d_2$ and $m.$  We obtain this result from 
 the above classification of $(A_1,A_2)$-invariant curves and  the 
following result of independent interest, which states roughly speaking that if 
a rational function $X$  is 
``a compositional left factor"  of some iterate of a rational function $A$, then $X$ is already a  factor  of $A^{\circ N}$, where $N$ is bounded in terms of degrees of $A$ and $X$.

\bt \l{5}  
There exists a function $\phi: \mathbb N\times \mathbb N\rightarrow \R$ with the following property. For any rational functions $A$ and $X$ such that the equality 
\be \l{decc} A^{\circ d}=X\circ R\ee holds for some rational function $R$ and $d\geq 1$, there exists 
$N\leq	 \phi(\deg A,\deg X)$ and a rational function $R'$ such that 
 \be \l{puzo} A^{\circ N}=X\circ R'\ee 
and $R=R'\circ A^{\circ (d-N)}$, if $d>N$. In particular, for any fixed rational function $A$ and  integer $n\geq 1$, up to the change $X\rightarrow X\circ \mu,$ where $\mu$ is a M\"obius transformation, 
there exist at most finitely many rational functions $X$ of degree $n$ such that \eqref{decc} holds for some rational function $R$ and $d\geq 1.$
\et

The paper is organized as follows. In the second and the third sections, we recall basic definitions and results related to orbifolds on Riemann surfaces, and
 review some of results of the papers \cite{semi} and \cite{lattes} describing the structure of solutions of functional equation \eqref{ssee} in rational functions. 
In the fourth section, we describe $(A_1,A_2)$-invariant, periodic, and preperiodic curves.   
In the fifth section, we prove results concerning the orbit density.

Finally, in the sixth section, we obtain quantitative versions of some results of the paper \cite{aol} concerning pairs of rational functions $A$ and $X$ such that for every $d\geq 1$ the algebraic curve $$A^{\circ d}(x)-X(y)=0$$ has a factor of genus zero or one. As an application, we prove Theorem \ref{5} and deduce from it the finiteness 
of the number of $(A_1,A_2)$-invariant curves of any given bi-degree $(d_1,d_2)$.

\vskip 0.1cm
\noindent {\it Acknowledgments}. The author would like to thank Dragos Ghioca, Laura DeMarco, Thomas Tucker, and Junyi Xie for helpful conversations.

\end{section}

\begin{section}{Orbifolds and generalized Latt\`es maps}

\begin{subsection}{Riemann surface orbifolds} 

{\it A Riemann surface  orbifold} is a pair $\f O=(R,\nu)$ consisting of a Riemann surface $R$ and a ramification function $\nu:R\rightarrow \mathbb N$, 
which takes the value $\nu(z)=1$ except at isolated points. 
For an orbifold $\f O=(R,\nu)$, 
 the {\it  Euler characteristic} of $\f O$ is the number
\be \l{euler} \chi(\f O)=\chi(R)+\sum_{z\in R}\left(\frac{1}{\nu(z)}-1\right),\ee
the set of {\it singular points} of $\f O$ is the set 
$$c(\f O)=\{z_1,z_2, \dots, z_s, \dots \}=\{z\in R \mid \nu(z)>1\},$$ and  the {\it signature} of $\f O$ is the set 
$$\nu(\f O)=\{\nu(z_1),\nu(z_2), \dots , \nu(z_s), \dots \}.$$ 
For   orbifolds $\f O_1=(R_1,\nu_1)$  and $\f O_2=(R_2,\nu_2)$, 
we write $ \f O_1\preceq \f O_2$ 
if $R_1=R_2$, and for any $z\in R_1$, the condition $\nu_1(z)\mid \nu_2(z)$ holds.

Let $\f O_1=(R_1,\nu_1)$  and $\f O_2=(R_2,\nu_2)$ be orbifolds
and let 
$f:\, R_1\rightarrow R_2$  be a holomorphic branched covering map. We say that $f:\,  \f O_1\rightarrow \f O_2$
is  a {\it covering map} 
{\it between orbifolds}
if for any $z\in R_1$ the equality 
\be \l{us} \nu_{2}(f(z))=\nu_{1}(z)\deg_zf\ee holds, where $\deg_zf$ is the local degree of $f$ at the point $z$.
If for any $z\in R_1$  
the weaker condition 
\be \l{uuss} \nu_{2}(f(z))\mid \nu_{1}(z)\deg_zf\ee
is satisfied  instead of  \eqref{us},   we say that $f:\,  \f O_1\rightarrow \f O_2$ 
is a {\it holomorphic map}
 {\it between orbifolds}.

A universal covering of an orbifold ${\f O}$
is a covering map between orbifolds \linebreak  $\theta_{\f O}:\,
\tt {\f O}\rightarrow \f O$ such that $\tt R$ is simply connected and $\tt {\f O}$ is non-ramified, that is, $\tt \nu(z)\equiv 1.$ 
If $\theta_{\f O}$ is such a map, then 
there exists a group $\Gamma_{\f O}$ of conformal automorphisms of $\tt R$ such that the equality 
$\theta_{\f O}(z_1)=\theta_{\f O}(z_2)$ holds for $z_1,z_2\in \tt R$ if and only if $z_1=\sigma(z_2)$ for some $\sigma\in \Gamma_{\f O}.$ 
A universal covering exists and 
is unique up to a conformal isomorphism of $\tt R$ whenever 
$\f O$ is {\it good}, that is,  distinct from the Riemann sphere with one ramified point or with two ramified points $z_1,$ $z_2$ such that $\nu(z_1)\neq \nu(z_2)$.   
 Furthermore, 
$\tt R$ is the unit disk $\mathbb D$ if and only if $\chi(\f O)<0,$ $\tt R$ is the complex plane $\C$ if and only if $\chi(\f O)=0,$ and $\tt R$ is the Riemann sphere $\P^1(\C)$ if and only if $\chi(\f O)>0$ (see e.g. \cite{fk}, Section IV.9.12).
Below we  always assume that considered orbifolds are good.
Abusing  notation, we  use the symbol $\tt {\f O}$ both for the
orbifold and for the  Riemann surface  $\tt R$.

Covering maps between orbifolds lift to isomorphisms between their univer\-sal coverings.
More generally, for any holomorphic map between orbifolds $f:\,  \f O_1\rightarrow \f O_2$ there exist 
a holomorphic map $F:\, \tt {\f O_1} \rightarrow \tt {\f O_2}$ and 
a homomorphism $\phi:\, \Gamma_{\f O_1}\rightarrow \Gamma_{\f O_2}$ such that the diagram 
\be \l{dia2}
\begin{CD}
\tt {\f O_1} @>F>> \tt {\f O_2}\\
@VV\theta_{\f O_1}V @VV\theta_{\f O_2}V\\ 
\f O_1 @>f >> \f O_2\ 
\end{CD}
\ee
commutes and 
for any $\sigma\in \Gamma_{\f O_1}$ the equality
\be \l{homm}  F\circ\sigma=\phi(\sigma)\circ F \ee holds.  
The holomorphic map $F$ is an isomorphism if and only if $f$ is a covering map between orbifolds (see \cite{semi}, Propo\-sition 3.1).

If $f:\,  \f O_1\rightarrow \f O_2$ is a covering map between orbifolds with compact supports, then  the Riemann-Hurwitz 
formula implies that 
\be \l{rhor} \chi(\f O_1)=d \chi(\f O_2), \ee
where $d=\deg f$. 
More generally,  if $f:\,  \f O_1\rightarrow \f O_2$ is a holomorphic map, then \be \l{iioopp} \chi(\f O_1)\leq \chi(\f O_2)\,\deg f, \ee and the equality is attained if and only if $f:\, \f O_1\rightarrow \f O_2$ is a covering map between orbifolds
(see \cite{semi}, Proposition 3.2).

Let $R_1$, $R_2$ be Riemann surfaces and 
$f:\, R_1\rightarrow R_2$ a holomorphic branched covering map. Assume that $R_2$ is provided with a ramification function $\nu_2$. In order to define a ramification function $\nu_1$ on $R_1$ so that $f$ would be a holomorphic map between orbifolds $\f O_1=(R_1,\nu_1)$ and $\f O_2=(R_2,\nu_2)$ 
we must satisfy condition \eqref{uuss}, and it is easy to see that
for any  $z\in R_1$ a minimal possible value for $\nu_1(z)$ is defined by 
the equality 
\be \l{rys} \nu_{2}(f(z))=\nu_ {1}(z)\GCD(\deg_zf, \nu_{2}(f(z)).\ee 
In case \eqref{rys} is satisfied for  any $z\in R_1$, we 
say that $f$ is {\it a  minimal holomorphic  map} 
between orbifolds 
$\f O_1=(R_1,\nu_1)$ and $\f O_2=(R_2,\nu_2)$.
It follows from the definition that for any orbifold $\f O=(R,\nu)$ and a holomorphic branched covering map \linebreak $f:\, R^{\prime} \rightarrow R$ there exists a unique orbifold structure  $\f O^{\prime}=(R^{\prime},\nu^{\prime})$  such that 
$f:\f O^{\prime}\rightarrow \f O$ is a minimal holomorphic map between orbifolds. 
We will denote the corresponding orbifold by $f^*\f O.$ Notice that any covering map between orbifolds $f:\,  \f O_1\rightarrow \f O_2$ is a  minimal holomorphic map. 

Minimal holomorphic maps between orbifolds possess the following fundamental property with respect to the operation of  composition (see \cite{semi}, Theorem 4.1).

\bt \l{serrr} Let $f:\, R^{\prime\prime} \rightarrow R^{\prime}$ and $g:\, R^{\prime} \rightarrow R$ be holomorphic branched covering maps, and  $\f O=(R,\nu)$ an orbifold. 
Then 

$$(g\circ f)^*\f O= f^*(g^*\f O).\eqno{\Box}$$
\et

Theorem \ref{serrr} implies  the following two corollaries (see   \cite{semi}, Corol\-lary 4.1 and Corollary 4.2).

\bc \l{serka0} Let $f:\, \f O_1\rightarrow \f O^{\prime}$ and $g:\, \f O^{\prime}\rightarrow \f O_2$ be minimal holomorphic maps (resp. covering maps) between orbifolds.
Then  $g\circ f:\, \f O_1\rightarrow \f O_2$ is  a minimal holomorphic map (resp. covering map). \qed
\ec

\bc \l{indu2}  Let $f:\, R_1 \rightarrow R^{\prime}$ and $g:\, R^{\prime} \rightarrow R_2$ be holomorphic branched covering maps, and  $\f O_1=(R_1,\nu_1)$ and  $\f O_2=(R_2,\nu_2)$
orbifolds. Assume that \linebreak $g\circ f:\, \f O_1\rightarrow \f O_2$ is  a minimal holomorphic map (resp. a co\-vering map). Then $g:\, g^*\f O_2\rightarrow \f O_2$  and  $f:\, \f O_1\rightarrow g^*\f O_2 $ are minimal holomorphic maps (resp. covering maps). \qed
\ec 
 
Most of orbifolds considered in this paper are defined on $\P^1(\C)$. 
For such orbifolds, we omit  the Riemann surface $R$ in the definition of $\f O=(R,\nu)$,
meaning that $R=\P^1(\C).$ 
Signatures of orbifolds on $\P^1(\C)$ with non-negative Euler characteristics and corresponding $\Gamma_{\f O}$ and $\theta_{\f O}$ can be described explicitly as follows. If $\f O$ is an orbifold distinct from the non-ramified sphere, then 
$\chi(\f O)=0$  if and only if the signature of $\f O$ 
belongs to the list
\be \l{list}\{2,2,2,2\} \ \ \ \{3,3,3\}, \ \ \  \{2,4,4\}, \ \ \  \{2,3,6\}, \ee and $\chi(\f O)>0$  if and only if
 the signature of $\f O$  belongs to the  list 
 \be \l{list2} \{l,l\}, \ \ l\geq 2,  \ \ \ \{2,2,l\}, \ \ l\geq 2,  \ \ \ \{2,3,3\}, \ \ \ \{2,3,4\}, \ \ \ \{2,3,5\}.\ee
Groups $\Gamma_{\f O}\subset Aut(\C)$ corresponding to orbifolds $\f O$ with signatures \eqref{list}  
are generated by translations of $\C$ by elements of some lattice $L\subset \C$ of rank two and the rotation $z\rightarrow  \v z,$ where $\v$ is an $n$th root of unity with $n$ equal to 2,3,4, or 6, such that  $\v L=L$ (see  \cite{mil2}, or \cite{fk}, 
Section IV.9.5).  
 Accordingly, the functions $\theta_{\f O}$ 
may be written in terms of the  corresponding
Weierstrass functions as $\wp(z),$ $\wp^{\prime }(z),$ $\wp^2(z),$  and $\wp^{\prime 2}(z).$  
Groups $\Gamma_{\f O}\subset Aut(\P^1(\C))$ corresponding to   orbifolds $\f O$ with signatures \eqref{list2} are the well-known finite subgroups 
 $C_l,$  $D_{2l},$  $A_4,$ $S_4,$ $A_5$ of $Aut(\P^1(\C))$, and the functions $\theta_{\f O}$ are Galois coverings of $\P^1(\C)$ by $\P^1(\C)$ of degrees 
$l$, $2l,$ $12,$ $24,$  $60,$ calculated for the first time by Klein in \cite{klein}.

\end{subsection}

\begin{subsection}{Functional equations and orbifolds}

 With each holomorphic map \linebreak $f:\, R_1\rightarrow R_2$ between compact Riemann surfaces, 
one can associate two orbifolds $\f O_1^f=(R_1,\nu_1^f)$ and 
$\f O_2^f=(R_2,\nu_2^f)$, setting $\nu_2^f(z)$  
equal to the least common multiple of local degrees of $f$ at the points 
of the preimage $f^{-1}\{z\}$, and $$\nu_1^f(z)=\frac{\nu_2^f(f(z))}{\deg_zf}.$$ By construction, 
 $$f:\, \f O_1^f\rightarrow \f O_2^f$$ 
is a covering map between orbifolds.
It is easy to  see that the covering map \linebreak $f:\, \f O_1^f\rightarrow \f O_2^f$ is minimal in the following sense. For any covering map between orbifolds $f:\, \f O_1\rightarrow \f O_2$ we have:
\be \l{elki+} \f O_1^f\preceq \f O_1, \ \ \ \f O_2^f\preceq \f O_2.\ee
Notice that for any orbifold $\f O$ the orbifolds $\f O_1^{\theta_{\f O}}$ and $\f O_2^{\theta_{\f O}}$ obviously are well defined even if $\t {\f O}$ is non-compact and satisfy 
\be \l{seq} \f O_1^{\theta_{\f O}}=\t {\f O}\ \ \ \ \ \ \ \f O_2^{\theta_{\f O}}=\f O.\ee

The orbifolds  defined above are useful for the study of the functional equation
\be \l{m} f\circ p=g\circ q, \ee 
where  
$$p:\, R\rightarrow C_1,  \ \ \ \ f:\, C_1\rightarrow \P^1(\C),\ \ \ \ q:\, R\rightarrow C_2, \ \ \ \ g:\, C_2\rightarrow \P^1(\C)$$ 
are holomorphic maps between compact Riemann surfaces.
We say that a solution $f,p,g,q$ of \eqref{m} is {\it good} if the fiber product of $f$ and $g$ has a unique component, 
and  $p:\, R\rightarrow C_1$ and $q:\, R\rightarrow C_2$ {\it have no non-trivial common compositional  right factor} in the following sense: the equalities 
$$ p= \tt p\circ  w, \ \ \ q= \tt q\circ w,$$ where $w:\, R \rightarrow {\tt R}$, $\tt p:\, {\tt R}\rightarrow C_1$, $\tt q:\,  {\tt R}\rightarrow C_2$ are holomorphic maps between compact Riemann surfaces, imply that $\deg w=1.$ 
In this notation, the following statement holds (see \cite{semi}, Theorem 4.2).

\bt \l{goodt} Let $f,p,g,q$ be a good solution of \eqref{m}.
Then the commutative diagram 
\be 
\begin{CD}
\f O_1^q @>p>> \f O_1^f\\
@VV q V @VV f V\\ 
\f O_2^q @>g >> \f O_2^f\ 
\end{CD}
\ee
consists of minimal holomorphic  maps between orbifolds. \qed  
\et

Good solutions admit the following characterization (see \cite{semi}, Lemma 2.1).

\bl \l{good} A 
solution $f,p,g,q$ of \e{m} is good whenever 
any two of the following three conditions are satisfied:

\begin{itemize}
\item the fiber product of $f$ and $g$ has a unique component,
\item $p$ and $q$ have no non-trivial common compositional right factor,
\item $ \deg f=\deg q, \ \ \ \deg g=\deg p.$  \qed
\end{itemize}
\el

Note that if $f$ and $g$ are rational functions, then
the fiber product of $f$ and $g$ has a unique component if and only if the algebraic
curve $f(x) - g(y) = 0$
is irreducible.

Finally, the following result (see \cite{aol}, Corollary 2.9 or \cite{fin}, Theorem 2.18) states  that ``gluing together'' two commutative diagrams corresponding to good solutions of \eqref{m} we obtain again a good solution of \eqref{m} (see the diagram below).
\be \l{gae}
\begin{CD}
\P^1(\C) @>B>> \P^1(\C)  @>W>> \P^1(\C) \\
@VV C V @VV D V @VV {V} V \\ 
\P^1(\C) @>A >> \P^1(\C)  @>U >> \P^1(\C)\,. 
\end{CD}
\ee

\vskip 0.2cm

\bt \l{sum+} Assume that the quadruples of rational functions 
$A,C,D,B$ and $U,D,V,W$ are  good solutions of \eqref{m}. 
Then the quadruple $U\circ A$, $C$, $V,$ $W\circ B$ is also  a good solution of \eqref{m}. \qed
\et

\end{subsection}

\begin{subsection}{Generalized Latt\`es maps}\l{gelatt}
We recall that a {\it Latt\`es map} $A$ is a rational function of degree at least two such that there exist a lattice $\Lambda$ of rank two in $\C$, an affine map $L=az+b$ on $\C$, 
 and a holomorphic map $\Theta:\C/\Lambda\rightarrow \P^1(\C)$,  such that 
$L(\Lambda)\subseteq \Lambda$ and the diagram 
\be \l{lla}
\begin{CD}
 \C/\Lambda @>az+b>>\C/\Lambda \\
@VV \Theta V @VV \Theta V\\ 
\P^1(\C) @> A >>\P^1(\C) \ 
\end{CD}
\ee
commutes (abusing the notation we will continue using the notation $az+b$ for the map on $\C/\Lambda$ induced by the affine map $az+b$ on  $\C$). 
 Equivalently, a Latt\`es map can be defined as a rational function $A$ of degree at least two such that $A:\f O\rightarrow \f O$ is a  covering self-map
for some orbifold $\f O$ (see \cite{mil2}). Thus, $A$ is a Latt\`es map  if there exists an orbifold $\f O$
such that  for any $z\in \P^1(\C)$ the equality 
\be \l{uu0} \nu(A(z))=\nu(z)\deg_zA\ee holds. By  formula \eqref{rhor}, such $\f O$ necessarily satisfies $\chi(\f O)=0.$ Furthermore, 
for a given function $A$ there might be at most one  orbifold such that \eqref{uu0} holds (see \cite{mil2} and \cite{lattes}, Theorem 6.1).  

Following \cite{lattes}, we say that a rational function $A$ of degree at least two is 
a {\it generalized Latt\`es map} if there exists an orbifold $\f O$, 
 distinct from the non-ramified sphere, 
such that  $A:\f O\rightarrow \f O$ is a minimal holomorphic self-map between orbifolds; that is, for any $z\in \P^1(\C)$, the equality 
\be \l{uu} \nu(A(z))=\nu(z)\GCD(\deg_zA,\nu (A(z)))\ee holds.  By inequality \eqref{iioopp}, such $\f O$ satisfies $\chi(\f O)\geq 0$.
Since condition \eqref{uu0} implies condition \eqref{uu}, any ordinary Latt\`es map 
is a generalized Latt\`es map. 
Note that if $\f O$ is the non-ramified sphere, then condition \eqref{uu} trivially holds for any rational function $A$.

In general, for a given function $A$ there might be several orbifolds $\f O$ satisfying \eqref{uu}, and even infinitely many such orbifolds. 
 For example, it is easy to see that $z^{\pm n}:\f O\rightarrow \f O$ is a minimal holomorphic map for any $\f O$ defined  by  $$\nu(0)=m, \ \ \ \ \nu(\infty)=m, \ \ \ \GCD(n,m)=1,$$ while
 $\pm T_{n}:\f O\rightarrow \f O$ is a minimal holomorphic map  for any $\f O$ defined by  the conditions $$
\nu(-1)=\nu(1)=2, \ \ \ \nu(\infty)=m, \ \ \ \  \GCD(n,m)=1.$$
 Nevertheless, the following statement holds (see \cite{lattes}, Theorem 1.2).

\bt \l{uni} 
Let $A$ be a rational function  of degree at least two not conjugate to $z^{\pm d}$ or $\pm T_d.$ Then there exists an orbifold $\f O_0^A$ such that $A:\, \f O_0^A\rightarrow \f O_0^A$
is a minimal holomorphic map between orbifolds, and for any orbifold $\f O$ such that 
$A:\, \f O\rightarrow \f O$ is a minimal holomorphic map between orbifolds, the relation $\f O\preceq \f O_0^A$ holds. Furthermore, $\f O_0^{A^{\circ l}}=\f O_0^A$ for any 
$l\geq 1$.  \qed
\et

It is clear that generalized Latt\`es maps are exactly rational functions for which the orbifold $\f O_0^A$ is distinct from the non-ramified sphere, completed
by the functions  $z^{\pm d}$ and $\pm T_d$ for which the orbifold $\f O_0^A$ is not defined. Furthermore,  ordinary Latt\`es maps are exactly  rational functions for which  $\chi(\f O_0^A)= 0$  (see \cite{lattes}, Lemma 6.4).
% and if $A$ is a  Latt\`es map, then the minimal holomorphic map $A:\f O_0^A\rightarrow \f O_0^A$ is a  covering map by Proposition \ref{p1}. 
Notice also that since a rational function $A$ is conjugate to $z^{\pm d}$ or $\pm T_d$ if and only 
 if some iterate $A^{\circ l}$, $l\geq 1,$   is conjugate to $z^{\pm ld}$ or $\pm T_{ld}$
 (see e.g. \cite{lattes}, Lemma 6.3),  Theorem \ref{uni} implies that $A$ is a generalized Latt\`es map if and only if some iterate $A^{\circ l}$, $l\geq 1,$  is a generalized Latt\`es map. Finally, notice that for a given rational function 
$A$ the orbifold $\f O_0^A$ can be effectively calculated from the branch data of $A$ (see \cite{lattes}, Section 6).

We recall that a rational function $A$ is called {\it special} if it is either a Latt\`es map, or it is conjugate to $z^{\pm n}$ or $\pm T_n.$
If $A$ is a generalized Latt\`es map, which is not special, then $\chi(\f O_0^A)>0$,
and the corresponding 
diagram \eqref{dia2} takes the form 
\be \l{dia3}
\begin{CD}
\P^1(\C) @>F>> \P^1(\C)\\
@VV\theta_{\f O_0^A}V @VV\theta_{\f O_0^A}V\\ 
\P^1(\C) @>A >>\P^1(\C)\,.
\end{CD}
\ee
Moreover, for such $A$ the homomorphism $\phi$ in \eqref{homm}
is an {\it automorphism}.
More precisely, the following statement holds
(see \cite{semi}, Theorem 5.1).

\bt \l{las}  Let $A$ and $F$ be rational functions of degree at least two, and  $\f O$ an orbifold with $\chi(\f O)>0$ such that 
$A:\, \f O \rightarrow \f O$ 
is a holomorphic map between orbifolds  and  the diagram 
\be \l{diaa}
\begin{CD}
\P^1(\C) @>F>> \P^1(\C)\\
@VV\theta_{\f O}V @VV\theta_{\f O}V\\ 
\f O @>A >> \f O\ 
\end{CD}
\ee
commutes.
Then the following conditions are equivalent:

\begin{enumerate}
\item The holomorphic map $A$ is a minimal holomorphic  map. 
\item  The homomorphism $\phi:\, \Gamma_{\f O}\rightarrow \Gamma_{\f O}$ defined by the equality 
\be \l{homoo} F\circ\sigma=\phi(\sigma)\circ F, \ \ \ \sigma\in \Gamma_{\f O},\ee is an automorphism of $\Gamma_{\f O}$.
\item The functions $\theta_{\f O}$, $F,$ $A,$ $\theta_{\f O}$  form a good solution of equation \eqref{m}. \qed
\end{enumerate}

\et

Finally, we need the following simple result  (see Lemma 6.6 of \cite{lattes}) imposing restrictions on ramification of generalized Latt\`es maps, and, more generally, on ramification of holomorphic coverings maps between orbifolds of positive Euler characteristic. 

\bl \l{toch} Let $A$ be a rational function of degree at least five, and $\f O_1$, $\f O_2$ orbifolds distinct from the non-ramified sphere such that  $A:\f O_1\rightarrow \f  O_2$ is a minimal holomorphic map between orbifolds. Assume that $\chi(\f O_1)\geq 0$. Then    $c(\f O_2)\subseteq c(\f O_2^A)$. 
\el

\end{subsection}

\end{section}

\begin{section}{Semiconjugate rational functions}
\begin{subsection}{Primitive solutions} 
Let $A$ and $B$ be rational functions of degree at least two. We recall that 
$B$ is said to be semiconjugate to $A$ if there exists a non-constant
rational function $X$ such that the equality
\be \l{i1}
A\circ X=X\circ B
\ee
holds. If $\deg X=1$, then $A$ and $B$ are conjugate in the usual sense.
We say that a solution $A,X,B$ of functional equation \eqref{i1} 
is {\it primitive} if $\C(B,X)=\C(x)$.  
By Lemma \ref{good}, a solution $A,X,B$ of  \eqref{i1}  is primitive if and only if the quadruple 
$$f=A, \ \ \ p=X, \ \ \  g=X, \ \ \ q=B$$  is  a good solution of \eqref{m}. Primitive solution are described as follows (see \cite{semi}, Theorem 6.1, or \cite{dyna}).

\bt \l{se1} Let  $A,X,B$ be a primitive solution of \e{i1} with $\deg X>1$. Then $\chi(\f O_1^X)\geq 0$, $\chi(\f O_2^X)\geq 0$, and 
the commutative diagram 
\be 
\begin{CD} \l{gooopa}
\f O_1^X @>B>> \f O_1^X\\
@VV X V @VV X V\\ 
\f O_2^X @>A >> \f O_2^X\ 
\end{CD}
\ee
consists of minimal holomorphic  maps between orbifolds. \qed
\et

In particular, Theorem \ref{se1} implies  that  if  $A,X,B$ is a primitive solution of \e{i1} with $\deg X>1$, then 
$A$ is necessarily a generalized Latt\`es map, and $X$ satisfies the condition $\chi(\f O_2^X)\geq 0,$
implying strong restrictions on $X$ (see \cite{gen}).
\vskip 0.2cm

\end{subsection}

\begin{subsection}{Elementary transformations}
Let $A$ be a rational function. 
For any decomposition $A=V\circ U,$ where $U$ and $V$ are rational functions, the 
rational function $\t A=U\circ V$ is called an {\it elementary transformation} of $A$, and rational functions $A$ and $B$ are called  {\it equivalent}  if there exists 
a chain of elementary transformations between $A$ and $B$. For a rational function $A$, we denote its equivalence class by $[A].$
Since for any M\"obius transformation $W$ the equality $A=(A\circ W)\circ W^{-1}$ holds, 
each equivalence class $[A]$ is a union of conjugacy classes. 
Moreover, an equivalence class $[F]$ 
contains infinitely many conjugacy classes if and only if 
$F$ is a flexible Latt\`es map (\cite{rec}). If $A$ is  a 
generalized Latt\`es map, then any elementary transformation of $A$ is  a 
generalized Latt\`es map (see \cite{lattes}, Theorem 4.1), implying that any $B\sim A$ is a  
generalized Latt\`es map.

The  connection between the relation $\sim$ and semiconjugacy is straightforward. Namely, for $\t A$ and $A$ as above 
the diagrams
$$ 
\begin{CD}
 \P^1(\C) @> A >> \P^1(\C) \\
 @VV U V @VV U  V \\
 \P^1(\C) @>  \t A >> \P^1(\C)\,,
\end{CD} \ \ \ \ \ \ \ \ \ \ \begin{CD}
 \P^1(\C) @>\t A >> \P^1(\C) \\
 @VV {V}  V @VV {V}  V \\
 \P^1(\C) @> A >> \P^1(\C)
\end{CD}
$$
commute, implying inductively that if
$A\sim \t A$, then $A$ is semiconjugate to $\t A$, and   $\t A$ is semiconjugate to $A.$ Moreover, the following statement, obtained by a direct calculation, is true (see \cite{lattes}, Lemma 3.1).

\bl \l{lem1} Let \be \l{chh} A\rightarrow A_1 \rightarrow A_2  \rightarrow \dots \rightarrow A_s\ee 
be a chain of elementary transformations, and   $U_i,$ $V_i,$ $1\leq i \leq s,$  rational functions such that 
$$A=V_1\circ U_1, \ \ \  A_i= U_i\circ V_i, \ \ \ \ \ 1\leq i\leq s,$$ 
and
\be \l{seqs}  U_{i}\circ V_{i}=V_{i+1}\circ U_{i+1},\ \ \ 1 \leq i \leq s-1.\ee
Then  the functions
\be \l{aga} U=U_s\circ U_{s-1}\circ \dots \circ U_{1}, \ \ \ \ V=V_{1}\circ \dots \circ V_{s-1}\circ V_s\ee
make the diagram 
\be 
\begin{CD} 
\P^1(\C) @> A>>\P^1(\C) \\ 
@V U  VV @VV U   V\\ 
 \P^1(\C) @> A_s>> \P^1(\C)
 \\ 
@V {V}  VV @VV {V}   V\\ 
 \P^1(\C) @> A>> \P^1(\C) 
\end{CD} 
\ee
commutative 
and 
satisfy the equalities 
$$V\circ U=A^{\circ s}, \ \ \ \ \  \ U\circ V=A_s^{\circ s}.\eqno{\Box}$$ 
\el

Non-primitive solutions of \eqref{i1} reduce to primitive ones by chains of elementary transformations 
(see \cite{semi} and \cite{lattes} for more detail). Below we only need the following statement.

\bp \l{eg} If $A,X,B$ is a solution of \e{i1} and $A$  is not a generalized Latt\`es map, then $B\sim A$ and  there exists a rational function $Y$ such that 
the diagram 
\be 
\begin{CD} 
\P^1(\C) @> B>>\P^1(\C) \\ 
@V X  VV @VV X   V\\ 
 \P^1(\C) @> A>> \P^1(\C)
 \\ 
@V {Y}  VV @VV {Y}   V\\ 
 \P^1(\C) @> B>> \P^1(\C), 
\end{CD} 
\ee
commutes, and the equalities  
\be \l{z} Y\circ X=B^{\circ d}, \ \ \ \ \  \ X\circ Y=A^{\circ d}\ee
hold for some $d\geq 0$. 
\ep
\pr 
In case $\deg X=1,$ the conlcusion of the proposition holds for $Y=X^{-1}$ and $d=0$.  Assume now that $\deg X>1.$ 
Since  $A$  is not a generalized Latt\`es map, it follows from Theorem \ref{se1} that the triple  $A,X,B$ is not a primitive solution of \e{i1}. 
 Therefore, by the L\"uroth theorem, $\C(B,X)=\C(U_1)$ for some rational function $U_1$ with $\deg U_1>1$, and hence
\be \l{of} B=V_1 \circ U_1, \ \ \ X=X_1\circ U_1\ee 
for some rational functions $X_1,$ $V_1$. Since equality \eqref{i1} implies the equality 
$$A\circ X_1=X_1 \circ (U_1\circ V_1),$$ the triple $A, X_1 ,U_1\circ V_1$ is  also  a solution of \e{i1}.
Moreover, this new solution  again is not primitive by 
Theorem \ref{se1}, implying that there exist rational functions $X_2,$ $V_2,$ $U_2$ such that 
$$U_1 \circ V_1=V_2 \circ U_2, \ \ \ X_1=X_2\circ U_2,$$ and 
$$A\circ X_2=X_2 \circ (U_2\circ V_2).$$ 
Continuing in this way and taking into account that 
$$\deg X>\deg X_1>\deg X_2\ \dots\,,$$ we obtain a chain of elementary transformations  between $A$ and $B$ and the representation $$X=U_s\circ U_{s-1}\circ \dots \circ U_{1}$$ as in Lemma \ref{lem1}, so the proposition follows from this lemma. \qed

\end{subsection}

\end{section}

\begin{section}{Invariant, periodic, and preperiodic curves.} 

\begin{subsection}{Invariant curves and semiconjugacies} 
Let $A_1,A_2$ be rational functions. We denote by $(A_1,A_2):(\P^1(\C))^2\rightarrow (\P^1(\C))^2$ the map given by the formula \be \l{for} (z_1,z_2)\rightarrow (A(z_1), A(z_2)).\ee 
We say that 
an irreducible algebraic curve $\c C$ in $(\P^1(\C))^2$  is $(A_1,A_2)$-{\it invariant} 
if  $(A_1,A_2)(\c C)= \c C,$ and $(A_1,A_2)$-{\it periodic} if  $$(A_1,A_2)^{\circ n}(\c C)=\c C$$
for some $n\geq 1.$ Finally, we say that $\c E$ is 
$(A_1,A_2)$-{\it preperiodic} if $(A_1,A_2)^{\circ l}(\c C)$ is periodic  
for some $ l\geq 1.$

The simplest $(A_1,A_2)$-invariant curves are vertical lines $x=a$, where $a$ is a fixed point of $A_1$, and 
horizontal lines $y=b$, where $b$ is a fixed point of $A_2$. Other invariant curves are described as follows.

\bt \l{in1}

Let $A_1,A_2$ be rational functions of degree at least two, and $\c C$ an irreducible $(A_1,A_2)$-invariant curve that is not a vertical or horizontal line.
Then the desingularization $\t{\c C}$ of $\c C$ has genus zero or one, and 
 there exist non-constant rational maps $X_1,X_2:\t{\c C}\rightarrow \P^1(\C)$ and   $B:\t{\c C}\rightarrow \t{\c C}$ such that 
 the diagram 
\be \l{du}
\begin{CD} 
(\t{\c C})^2 @>(B,B)>> (\t{\c C})^2 \\ 
@V (X_1,X_2)  VV @VV  (X_1,X_2) V\\ 
 (\P^1(\C))^2 @>(A_1,A_2)>> (\P^1(\C))^2 
\end{CD} 
\ee
commutes and the map $t\rightarrow (X_1(t),X_2(t))$ is a generically one-to-one parametrization of $\c C.$   Finally, unless both $A_1,$ $A_2$ are Latt\`es maps, $\t{\c C}$ has genus zero.
\et

\pr Let $\t{\c C}$ be  the desingularization of $\c C$, and  
 $\pi: \t{\c C}\rightarrow \c C$ the desingularization map. We set 
$$X_1=x\circ \pi, \ \ \ \ \ X_2=y\circ \pi,$$ where $x,y:(\P^1(\C))^2\rightarrow \P^1(\C)$ are the projections on the first and on the second coordinate correspondingly. 
Since the map $(X_1,X_2): \t{\c C}\rightarrow \c C$ is an isomorphism off a finite set of points, 
the  map $(A_1,A_2): \c C \rightarrow \c C$  
 lifts to a rational map $B:\t{\c C}\rightarrow \t{\c C}$ which makes diagram \eqref{du} commutative.  
Furthermore, since $\c C$ is not a vertical or horizontal line,  
$X_1$ and $X_2$ are non-constant, implying by \eqref{du} that  
 $$\deg A_1=\deg A_2=\deg B.$$ In particular, $\deg B\geq 2$. 
It follows now from the Riemann-Hurwitz formula
$$2g(\t{\c C})-2=(2g(\t{\c C})-2)\deg B+\sum_{P\in \t{\c C}}(e_p-1)$$ that $g(\t{\c C})\leq 1.$  
Finally, if $g(\t{\c C})=1$, then  $A_1$ and $A_2$ are Latt\`es maps. Indeed, in this case 
$\t{\c C}=\C/\Lambda$ for some lattice $\Lambda$, and $B:\C/\Lambda\rightarrow \C/\Lambda$ is induced by an affine map. Thus, diagram \eqref{du} consists of a pair of 
 diagrams of the form \eqref{lla}. 
 \qed

\br
Note that Theorem \ref{in1} implies in particular that if $\deg A_1\neq \deg A_2$, then any  $(A_1,A_2)$-invariant curve
is  a vertical or horizontal line.
\er

\br \l{br1} 
For  an arbitrary field of characteristic zero $K$ and rational functions $A_1$ and $A_2$ defined over $K$   the notion of invariant curve is defined in the same way as above. Furthermore, it is easy to see that if $K$ is 
algebraically closed, then an analogue of Theorem \ref{in1} remains true over $K$. 

Below we will consider only fields $K$ that are subfields of $\C.$ 
For the problems  considered in this paper, such a restriction does not lead to the loss of generality since 
$A_1$ and $A_2$ defined over a field of characteristic zero actually are defined over a finitely generated extension of $\Q$, and such an extension can be embedded into $\C.$
 Taking this into account, we do not consider the case $K\neq \C$ separately till the fifth section. Note that 
the assumption $K\subset \C$ allows us in particular to continue using the notion  of  a generalized Latt\`es map.  

\er

The following lemma relates periodic curves for pairs of semiconjugate maps.

\bl \l{lk} Let $A_1,$ $A_2,$ $B_1,$ $B_2,$ $X_1,$ $X_2$ be non-constant rational functions such that 
the  diagram 
\be\l{kot}
\begin{CD} 
(\P^1(\C))^2 @>(B_1,B_2)>>(\P^1(\C))^2 \\ 
@V (X_1,X_2)  VV @VV  (X_1,X_2) V\\ 
 (\P^1(\C))^2 @>(A_1,A_2)>> (\P^1(\C))^2 
\end{CD} 
\ee
commutes. Then 
for any irreducible $(A_1,A_2)$-periodic (resp. preperiodic) curve $\c C$ there exists  
 an irreducible $(B_1,B_2)$-periodic  (resp. preperiodic) curve $\c C'$ such that  $\c C=(X_1,X_2)(\c C')$. 
\el
\pr 
For any  irreducible curve  $\c C$ in $(\P^1(\C))^2$ the preimage 
 $\c E=(X_1,X_2)^{-1}(\c C)$ is a union of irreducible curves, and any irreducible component ${\c C}'$ of $\c E$ satisfies $ (X_1,X_2)({\c C}')=\c C.$ Furthermore, if $\c C$ satisfies $(A_1,A_2)^{\circ n}(\c C)=\c C$,  then  
 $\c E$ satisfies $(B_1,B_2)^{\circ n}(\c E)\subseteq \c E,$
implying that all components of $\c E$ are $(B_1,B_2)$-preperiodic and at least one of these components is  $(B_1,B_2)$-periodic.
 Similarly, if 
 $\c C$ is $(A_1,A_2)$-preperiodic, then any component ${\c C}'$ of $\c E$ is  $(B_1,B_2)$-preperiodic.
\qed

\vskip 0.2cm

Assuming that at least one $(A_1,A_2)$-invariant curve $\c C$ is known, Theorem \ref{in1} combined with Lemma \ref{lk} permits to reduce describing $(A_1,A_2)$-periodic curves for a pair of functions $A_1,$ $A_2$ to describing 
$(B,B)$-periodic curves for a single  function $B$.

\bc \l{cor1} Let $A_1,A_2$ be rational functions of degree at least two that are not Latt\`es maps, and $\c B$    
a fixed irreducible $(A_1,A_2)$-invariant curve that is not a vertical or horizontal line. Then there exist rational functions $X_1,X_2,B$ 
such that diagram \eqref{ddii} commutes, the map $t\rightarrow (X_1(t),X_2(t))$ is a parametrization of $\c B,$ 
and  any irreducible $(A_1,A_2)$-periodic (resp. preperiodic) curve ${\c C}$ is the   $(X_1,X_2)$-image of some
 irreducible $(B,B)$-periodic  (resp. preperiodic) curve $\c C'$.
\qed
\ec

\end{subsection}

\begin{subsection}{The case where $A_1$, $A_2$ are not generalized Latt\`es maps} 
In this section, we describe $(A_1,A_2)$-invariant, periodic, and preperiodic curves in the case where $A_1$, $A_2$ are not generalized Latt\`es maps. 
%We start by proving Theorem \ref{1} from the introduction. 
\vskip 0.2cm
\noindent{\it Proof of Theorem \ref{1}.}
It was already mentioned in the introduction, that for any rational functions $X_1,X_2,A,B$ that make diagram \eqref{ddii} commutative, the map $t\rightarrow (X_1(t),X_2(t))$ is a parametrization of some $(A_1,A_2)$-invariant curve $\c C.$

In the other direction, assume that $\c C$ is  an $(A_1,A_2)$-invariant curve. 
Then by Theorem \ref{in1} there exist rational  functions $X_1,$ $X_2$, $B$  such that  diagram \eqref{ddii} commutes and the map $t\rightarrow (X_1(t),X_2(t))$ is a parametrization of $\c C.$   Furthermore, since $A_1$ and $A_2$ are not generalized Latt\`es maps, it follows from 
Proposition \ref{eg} that  there exist rational functions $Y_i$, $i=1,2$, 
  such that  
the diagram 
\be \l{gin} 
\begin{CD} 
(\P^1(\C))^2 @> (B,B) >>(\P^1(\C))^2 \\ 
@V (X_1,X_2)  VV @VV  (X_1,X_2)   V\\ 
 (\P^1(\C))^2 @> (A_1,A_2) >> (\P^1(\C))^2
 \\ 
@V (Y_1,Y_2)  VV @VV (Y_1,Y_2)    V\\ 
 (\P^1(\C))^2 @> (B,B)>> (\P^1(\C))^2 
\end{CD} 
\ee
commutes and the equalities 
\be \l{egi} X_i\circ Y_i=A_i^{\circ d_i}, \ \ \ \ \  \ Y_i\circ X_i=B^{\circ d_i}, \ \ \ \ i=1,2,\ee hold for some $d_1,d_2\geq 0.$ 

Let us show that modifying $Y_1$ and $Y_2$ we may assume that $d_1=d_2$. Suppose, say, that $d_2\geq d_1$. Setting $d=d_2$ and completing diagram \eqref{gin} to the diagram 
$$
\begin{CD} 
 (\P^1(\C))^2 @> (B,B)>> (\P^1(\C))^2
 \\ 
@V (X_1,X_2)  VV @VV (X_1,X_2)    V\\ 
 (\P^1(\C))^2 @> (A_1,A_2)>> (\P^1(\C))^2\\
@V (Y_1,Y_2)  VV @VV  (Y_1,Y_2)   V\\ 
 (\P^1(\C))^2 @> (B,B)>> (\P^1(\C))^2
 \\ 
@V (B^{\circ(d_2-d_1)},z)  VV @VV  (B^{\circ(d_2-d_1)},z)  V\\ 
 (\P^1(\C))^2 @> (B,B)>> (\P^1(\C))^2,
\end{CD} 
$$
we see that for the rational functions 
$$\t Y_1=B^{\circ(d_2-d_1)}\circ Y_1, \ \ \ \ \ \ \ \ \t Y_2= Y_2$$ 
diagram \eqref{gin} still commutes. Moreover, 
$$X_1\circ  \t Y_1=X_1\circ B^{\circ(d_2-d_1)}\circ Y_1=A_1^{\circ(d_2-d_1)}\circ X_1\circ Y_1=A_1^d,$$
$$X_2\circ  \t Y_2=X_2\circ Y_2=A_2^d,$$ and 
$$\t Y_i\circ X_i=B^{\circ d}, \ \ \ \ i=1,2.\eqno{\Box}$$

\vskip 0.2cm

%In applications, it is often desirable to know a description of $(A_1, A_2)$-periodic and preperiodic curves rather than invariant ones. In fact, the description of such
%curves is somewhat easier. Specifically, the following statement holds.

\bt \l{2} Let $A_1$, $A_2$ be rational functions of degree at least two that are not generalized Latt\`es maps, and 
$\c C$ an irreducible algebraic curve  in $(\P^1(\C))^2$ that is not a vertical or horizontal line. Then $\mathcal C$ is $(A_1,A_2)$-periodic  if and only if 
 there exist  rational functions  $X_1,$ $X_2,$ $Y_1,$ $Y_2$ such that
the equalities 
\be \l{en1+} X_1\circ Y_1=A_1^{\circ d}, \ \ \ \ \ \ \ 
 X_2\circ Y_2=A_2^{\circ d},\ee
 \be\l{en2+} Y_1\circ X_1=Y_2\circ X_2\ee hold for some $d\geq 0$, and   
the map $t\rightarrow (X_1(t),X_2(t))$ is a parametrization of  $\c C$.  On the other hand, $\c C$ is  $(A_1,A_2)$-preperiodic  if and only if 
 there exist  rational functions as above such that 
$\c C$ is a component of the curve 
 $Y_1(x)-Y_2(y)=0.$ 
\et

\noindent{\it Proof of Theorem \ref{2}.}
If $(A_1,A_2)^{\circ l}(\c C)=\c C$ for some $l\geq 1$, then by Theorem \ref{1} there exist rational functions $X_1,X_2,Y_1,Y_2,B$ such that the diagram 
\be  \l{xyxy}
\begin{CD} 
(\P^1(\C))^2 @>(B,B)>>(\P^1(\C))^2 \\ 
@V (X_1,X_2)  VV @VV  (X_1,X_2) V\\ 
 (\P^1(\C))^2 @>(A_1^{\circ l},A_2^{\circ l})>> (\P^1(\C))^2
\\ 
@V (Y_1,Y_2)  VV @VV  (Y_1,Y_2) V\\ 
 (\P^1(\C))^2 @>(B,B)>> (\P^1(\C))^2 
\end{CD}
\ee
commutes, the equalities 
 \be\l{ne1}  X_1\circ Y_1=A_1^{\circ d_0l}, \ \ \ \ \ \ \ 
 X_2\circ Y_2=A_2^{\circ d_0l},\ee
 \be\l{ne2} Y_1\circ X_1=Y_2\circ X_2=B^{\circ d_0}\ee hold for some $d_0\geq 0$, 
and $t\rightarrow (X_1(t),X_2(t))$ is a parametrization of $\c C.$ 
Thus, \eqref{en1+} and \eqref{en2+} hold for $d=ld_0.$

 On the other hand, if \eqref{en1+} and \eqref{en2+} hold, then setting \be \l{b} B= Y_1\circ X_1=Y_2\circ X_2\ee we see that the diagram 
\be \l{gopa}
\begin{CD} 
(\P^1(\C))^2 @> (B,B)>>(\P^1(\C))^2 \\ 
@V (X_1,X_2)  VV @VV  (X_1,X_2)   V\\ 
 (\P^1(\C))^2 @> (A_1^{\circ d},A_2^{\circ d})>> (\P^1(\C))^2
\end{CD}
\ee
commutes, implying that the curve $\f C$ parametrized by the map  $t\rightarrow (X_1(t),X_2(t))$ satisfies 
$(A_1,A_2)^{\circ d}(\c C)=\c C$.
This proves the first part of the theorem.

Assume now that $\c C'$ is an $(A_1,A_2)$-preperiodic curve. Then there exists 
a curve $\c C$ such that $(A_1,A_2)^{\circ l}(\c C)=\c C$ for some $l\geq 1$
and $\c C'$ is contained in the preimage of  $\c C$ under the map $(A_1,A_2)^{\circ s}$ for some $s\geq 0.$ Therefore, by the already proved part of the theorem, 
$\c C'$ is a component of the curve 
$$(Y_1\circ A_1^{\circ s})(x)- (Y_2\circ A_2^{\circ s})(y)=0$$ for some rational functions $Y_1,Y_2$ satisfying \eqref{xyxy}, \eqref{ne1}, \eqref{ne2}. Moreover, since the equality $$(A_1,A_2)^{\circ (l+s_0)}(\c C')=(A_1,A_2)^{\circ s_0}(\c C')$$ implies that
$$(A_1,A_2)^{\circ (l+s)}(\c C')=(A_1,A_2)^{\circ s}(\c C')$$ for any $s\geq s_0$, 
without loss of generality we may assume that $s=tl$ for some $t\geq 1$. Thus,   $\c C'$ is a component of the curve 
$$Y_1'(x)- Y_2'(y)=0,$$ where  
$$Y_1'=Y_1\circ A_1^{\circ tl}, \ \ \ \ Y_2'=Y_2\circ A_2^{\circ tl},$$ and  the functions $Y_1',$ $Y_2'$    
 satisfy the required conditions \eqref{en1+} and \eqref{en2+}, since 
\be \l{xxx1}  X_i\circ Y_i'=A_i^{\circ d_0l}\circ A_i^{\circ tl}=A_i^{\circ (d_0+t)l}, \ \ \ \ i=1,2,\ee
and 
\be \l{xxx2} Y_i'\circ X_i=Y_i\circ A_i^{\circ tl}\circ X_i=Y_i\circ X_i\circ B^{\circ t}=B^{\circ d_0}\circ B^{\circ t}=B^{\circ (d_0+t)}, \ \ \ \ i=1,2.\ee 

Lastly, if \eqref{en1+} and \eqref{en2+} hold, then for $B$ defined by formula \eqref{b} the diagram 
\be 
\begin{CD} 
(\P^1(\C))^2 @> (A_1^{\circ d},A_2^{\circ d})>>(\P^1(\C))^2 \\ 
@V (Y_1,Y_2)  VV @VV  (Y_1,Y_2)   V\\ 
 (\P^1(\C))^2 @> (B,B)>> (\P^1(\C))^2
\end{CD}
\ee
commutes. Therefore, curve \eqref{giv} satisfies $(A_1,A_2)^{\circ d}(\c E)\subseteq \c E,$
implying that every component of $\c E$ is preperiodic. \qed

\br \l{rem1}
Note that for every $(A_1,A_2)$-invariant curve $\c C$ we can 
find rational functions  $X_1,X_2,Y_1,Y_2,B$ satisfying conditions 1)-3) of Theorem  \ref{1} and the additional condition
that the parametrization $t\rightarrow (X_1(t),X_2(t))$ of $\c C$  is generically one-to-one, or equivalently that $\C(X_1,X_2)=\C(z).$
  Indeed, the functions $Y_1$ and $Y_2$ in the proof of the necessity are constructed from the functions $X_1$ and $X_2$ provided by Theorem \ref{in1}, and these functions  satisfy the required condition. 
 However, arbitrary rational functions satisfying \eqref{xx},  \eqref{en1}, \eqref{en2} 
do not necessarily satisfy condition $\C(X_1,X_2)=\C(z).$ A similar remark holds for Theorem \ref{2}. 
\er

\br \l{br2} Note that if under the assumptions of Theorem \ref{1} the functions $A_1$, $A_2$ are defined 
over an algebraically closed field $K\subset \C$, then we can assume that the functions $X_1,X_2,Y_1,Y_2,B$ are also defined over $K$. Indeed, for $X_1,X_2,B$ this is a corollary of Theorem \ref{in1} (see Remark \ref{br1}). On the other hand, if $X_1,X_2,B$ are defined over $K$, then $Y_1$, $Y_2$ are also defined over $K$ since their  coefficients are given by a linear system of equations  over $K$ obtained from the second group of equalities in \eqref{egi}.   
 A similar remark holds for Theorem \ref{2}. 

\er

\end{subsection}

\begin{subsection}{\l{ca} The case where $A_1=A_2$} In this section we provide an alternative description of $(A_1, A_2)$-invariant, periodic, and preperiodic  curves in
the special case $A_1 = A_2 = A$ in terms of functions commuting with $A$ or with some iterate of $A$.

\vskip 0.2cm

\noindent{\it Proof of Theorem \ref{1+}.} If $\c C$ is $(A,A)$-invariant, then  
applying Theorem \ref{1} we can  find  rational functions $X_1,X_2,Y_1,Y_2,B$ such that 
the diagram 
\be \l{as} 
\begin{CD} 
(\P^1(\C))^2 @> (B,B)>>(\P^1(\C))^2 \\ 
@V (X_1,X_2)  VV @VV  (X_1,X_2)   V\\ 
 (\P^1(\C))^2 @> (A,A)>> (\P^1(\C))^2
 \\ 
@V (Y_1,Y_2)  VV @VV (Y_1,Y_2)    V\\ 
 (\P^1(\C))^2 @> (B,B)>> (\P^1(\C))^2, 
\end{CD} 
\ee
commutes, 
the equalities 
\be \l{begi} X_i\circ Y_i=A^{\circ d_0}, \ \ \ \ \  \ Y_i\circ X_i=B^{\circ d_0}, \ \ \ \ i=1,2,\ee hold for some $d_0\geq 0,$ and  $t\rightarrow (X_1(t),X_2(t))$ is a parametrization of  $\c C$. 
Completing now diagram \eqref{as} to the diagram 
\be \l{qwer} 
\begin{CD} 
(\P^1(\C))^2 @> (A,A)>>(\P^1(\C))^2 \\ 
@V (Y_1,Y_1)  VV @VV  (Y_1,Y_1)   V\\ 
 (\P^1(\C))^2 @> (B,B)>> (\P^1(\C))^2
 \\ 
@V (X_1,X_2)  VV @VV (X_1,X_2)    V\\ 
 (\P^1(\C))^2 @> (A,A)>> (\P^1(\C))^2\\
@V (Y_1,Y_2)  VV @VV  (Y_1,Y_2)   V\\ 
 (\P^1(\C))^2 @> (B,B)>> (\P^1(\C))^2
 \\ 
@V (X_1,X_1)  VV @VV (X_1,X_1)    V\\ 
 (\P^1(\C))^2 @> (A,A)>> (\P^1(\C))^2,
\end{CD} 
\ee
and setting 
\be \l{ggpp} U_1=X_1\circ Y_1, \ \ \ \  U_2=X_2\circ Y_1, \ \ \ \ V_1=X_1\circ Y_1,  \ \ \  \ V_2=X_1\circ Y_2,\ee we see that the diagram 
\be  \l{xx1}
\begin{CD} 
(\P^1(\C))^2 @>(A,A)>>(\P^1(\C))^2 \\ 
@V (U_1,U_2)  VV @VV  (U_1,U_2) V\\ 
 (\P^1(\C))^2 @>(A,A)>> (\P^1(\C))^2
\\ 
@V ({V}_1,{V}_2)  VV @VV  ({V}_1,{V}_2) V\\ 
 (\P^1(\C))^2 @>(A,A)>> (\P^1(\C))^2 
\end{CD}
\ee
 commutes, implying that $U_1,$ $U_2,$ $V_1,$ $V_2$    commute with  $A$.

Furthermore, we 
have:
$$V_i\circ U_i=X_1\circ Y_i\circ X_i\circ Y_1=X_1\circ B^{\circ d_0}\circ Y_1=A^{\circ d_0}\circ X_1\circ Y_1=A^{\circ 2d_0},\ \ \ \ \ i=1,2,$$
and 
$$U_i\circ V_i=X_i\circ Y_1\circ X_1\circ Y_i=X_i\circ B^{\circ d_0}\circ Y_i=A^{\circ d_0}\circ X_i\circ Y_i=A^{\circ 2d_0},\ \ \ \ \ i=1,2,$$ implying that equalities \eqref{en11} and \eqref{en21} hold for $d=2d_0$.
Finally, since obviously $(Y_1,Y_1)(\Delta)=\Delta,$ the equality 
\be \l{ccuu} (U_1,U_2)(\Delta)=(X_1,X_2)(\Delta)=\c C\ee 
holds, that is,  $t\rightarrow (U_1(t),U_2(t))$ is a parametrization of  $\c C$. This proves the necessity. 

The sufficiency follows merely from the commutativity of the top square of diagram \eqref{xx1}, which in turn follows from the assumption that $U_1$ and $U_2$ commute with $A$. \qed

\br Note that for $A_2\neq A_1$  an analogue of diagram \eqref{qwer} is obtained by 
changing $(Y_1,Y_1)$ to 
$(Y_1, Y_2)$ and $(X_1,X_1)$ to 
$(X_1, X_2).$ Nevertheless, equality \eqref{ccuu} does not hold anymore since 
$(Y_1,Y_2)(\Delta)\neq \Delta.$

\er

\bt \l{2+} Let $A$ be a rational function of degree at least two that is not a generalized Latt\`es map, and 
$\c C$ an irreducible algebraic curve  in $(\P^1(\C))^2$  that is not a vertical or horizontal line. Then $\c C$ is $(A,A)$-periodic  if and only if 
 there exist rational functions $U_1,$ $U_2,$ $V_1,$ $V_2$ commuting with some iterate of $A$ such that 
the equalities 
\be \l{en12} U_1\circ V_1=U_2\circ V_2=A^{\circ d},\ee
 \be\l{en22} V_1\circ U_1=V_2\circ U_2=A^{\circ d}\ee hold for some $d\geq 0$  and   
the map $t\rightarrow (U_1(t),U_2(t))$ is a parametrization of  $\c C$.  On the other hand, $\c C$ is  $(A,A)$-preperiodic  if and only if 
 there exist  rational functions as above such that 
$\c C$ is a component of the curve 
 $V_1(x)-V_2(y)=0.$
\et
\pr The first part of the theorem follows directly from Theorem \ref{1+}, 
so we only must prove the second part.

The sufficiency follows from the commutativity of the diagram 
$$
\begin{CD} 
(\P^1(\C))^2 @> (A^{\circ d},A^{\circ d})>>(\P^1(\C))^2 \\ 
@V ({V}_1,{V}_2)  VV @VV  ({V}_1,{V}_2)  V\\ 
 (\P^1(\C))^2 @> (A^{\circ d},A^{\circ d})>> (\P^1(\C))^2
\end{CD}
$$
in the same way as in the proof of Theorem \ref{2}. 
To prove the necessity, let us observe that equality  \eqref{en21} implies, in the notation of the proof of Theorem \ref{1+}, that the invariant curve \eqref{ccuu} 
is a component of the curve 
defined by the equation 
$$V_1(x)-V_2(y)=0.$$  Therefore,  if $\c C'$ is an $(A,A)$-preperiodic curve,  then $\c C'$ is a component of the curve 
$$(V_1\circ A^{\circ s})(x)- (V_2\circ A^{\circ s})(y)=0$$ for some $s\geq 0$ and rational functions $V_1,V_2$, which commute with $A^{\circ l}$ for some 
 $l\geq 1$ and 
satisfy \eqref{en12}, \eqref{en22}. Furthermore, as 
in the proof of Theorem \ref{2}, without loss of generality we may assume that $\c C'$ is a component of $$V_1'(x)-V_2'(y)=0,$$ where 
$$V_1'=V_1\circ A^{\circ tl}, \ \ \ \ V_2'=V_2\circ A^{\circ tl}$$ fo some $t\geq 1.$ Finally, $V_1'$ and $V_2'$ commute with  $A^{\circ l}$
and satisfy 
$$ U_i\circ {V}_i'=A^{\circ (d+tl)},\ \ \ i=1,2,$$
$${V}_i'\circ U_i={V}_i\circ A^{\circ tl} \circ U_i=A^{\circ tl}\circ {V}_i\circ U_i=A^{\circ tl}\circ A^{\circ d}=A^{\circ (d+tl)}, \ \ \ i=1,2.\eqno{\Box}$$

\br \l{15}
Note that since the functions $U_1,U_2,V_1,V_2$ in Theorem \ref{1+} and Theorem \ref{2+} commute with some iterate of $A$, and $A$ is not special, it follows from the Ritt theorem about commuting rational functions (see \cite{r}) that each of the functions $U_1,U_2,V_1,V_2$ has a common iterate with $A$.
\er

\br \l{br3} Note that if under the assumptions of Theorem \ref{1+} the function $A$ is defined 
over an algebraically closed field $K\subset \C$, then we can assume that the functions $U_1,U_2,V_1,V_2$ are also defined over $K$. Indeed, the functions $U_1,U_2,V_1,V_2$ are given by equality \eqref{ggpp}, and the functions $X_1,X_2,Y_1,Y_2$ can be  defined over $K$ by Remark \ref{br2}. A similar remark holds for Theorem \ref{2+}. 

\er

\end{subsection}

\begin{subsection}{Description of $(A_1,A_2)$-invariant curves for non-special $A_1,$ $A_2$}

In this section, we show that describing $(A_1,A_2)$-periodic and preperiodic curves for non-special $A_1,$ $A_2$ can be reduced to the case 
where  $A_1$ and $A_2$  are not generalized Latt\`es maps.

\bl \l{ser} 
Let $U,V,X$  be rational functions such that $X=U\circ V$. Then $ \f O_2^U\preceq  \f O_2^X.$ Moreover, if  
$ \f O_2^U=  \f O_2^X$, then $ \f O_2^V\preceq \f O_1^U.$

\el
\pr Since $X:\, \f O_1^X\rightarrow \f O_2^X$ is a covering map, it follows from Corollary \ref{indu2} that 
\be \l{cove} U:\, U^*\f O_2^X\rightarrow \f O_2^X, \ \ \ \ V:\, \f O_1^X\rightarrow  U^*\f O_2^X\ee
are covering maps. Therefore,  since 
\be \l{comaps} U:\, \f O_1^U\rightarrow \f O_2^U, \ \ \ \ V:\, \f O_1^V\rightarrow  \f O_2^V\ee
also  are covering maps, the relation $ \f O_2^U\preceq  \f O_2^X$ holds by \eqref{elki+}.
Moreover,  in addition, we see that 
 \be \l{rela2} \f O_1^U\preceq  U^*\f O_2^X, \ \ \ \f O_2^V\preceq  U^*\f O_2^X.\ee

It follows from formula \eqref{rhor} applied to the first covering in \eqref{cove} that 
$$\chi(U^*\f O_2^X)=\deg U\cdot\chi(\f O_2^X).$$ Since, on the other hand,  
$$\chi(\f O_1^U)=\deg U\cdot\chi(\f O_2^U),$$ we see that if $ \f O_2^U=  \f O_2^X$, then \be \l{kva} \chi(\f O_1^U)=\chi(U^*\f O_2^X).\ee Since for any pair of orbifolds satisfying $\t{\f O}\preceq \f O$ the equality $\chi(\t{\f O})=\chi(\f O)$ 
holds if and only if $\t{\f O}=\f O$, equality \eqref{kva} and the first relation in \eqref{rela2} imply  that 
$\f O_1^U= U^*\f O_2^X$. It follows now from the second  relation in \eqref{rela2} that   $ \f O_2^V\preceq \f O_1^U$. \qed

%The following statement is ``the orbifold counterpart'' of Theorem \ref{3}. 

\bt \l{33} Let $A$ be a non-special rational function of degree at least two, and $B$ a rational function that makes the diagram 
\be \l{asa}
\begin{CD} 
\P^1(\C) @> B>>\P^1(\C) \\ 
@V \theta_{\f O_0^A}  VV @VV \theta_{\f O_0^A}   V\\ 
 \P^1(\C) @> A>> \P^1(\C). 
\end{CD} 
\ee
commutative. Then the orbifold $\f O_0^B$ is the non-ramified sphere. 
\et
\pr Let us complete diagram \eqref{asa} to the diagram 
\be \l{46} 
\begin{CD} 
\t{\f O}_0^B @> C>>\t{\f O}_0^B \\ 
@V \theta_{\f O_0^B}  VV @VV \theta_{\f O_0^B}   V\\ 
 \P^1(\C) @> B>> \P^1(\C)
 \\ 
@V \theta_{\f O_0^A}  VV @VV \theta_{\f O_0^A}   V\\ 
 \P^1(\C) @> A>> \P^1(\C)\,,
\end{CD} 
\ee 
and set \be \l{ru0}X=\theta_{\f O_0^A} \circ \theta_{\f O_0^B}.\ee 
Let us observe first that $\t{\f O}_0^B=\P^1(\C)$, implying that the functions $\theta_{\f O_0^B}$ and $X$ are rational. Indeed, since $\chi(\f O_0^B)\geq 0,$ 
otherwise $\t{\f O}_0^B =\C$, $C=az+b$ for some $a,b\in \C$, and
$\theta_{\f O_0^B}$ and $X$ are doubly periodic meromorphic function with respect to some lattice $\Lambda$. Therefore, in this case  diagram \eqref{lla} commutes for some holomorphic function $\Theta$,  in contradiction with the assumption that $A$ is not a Latt\`es map.

Since the quadruples $A, \theta_{\f O_0^A}, \theta_{\f O_0^A}, B$ and $B, \theta_{\f O_0^B}, \theta_{\f O_0^B}, C$
are good solutions of 
\eqref{m} by Theorem \ref{las}, the quadruple $A,X,X,C$ is also a good solution of \eqref{m} by Theorem \ref{sum+}, implying that $A:\f O_2^X\rightarrow \f O_2^X$ is a minimal holomorphic map by Theorem \ref{goodt}. 
Therefore, \be \l{kroll} \f O_2^X \preceq \f O_0^A\ee by Theorem \ref{uni}.
Since $$\f O_2^{ \theta_{\f O_0^A}}\preceq \f O_2^X $$ by the first part of Lemma \ref{ser}
and 
\be \l{plm} \f O_2^{ \theta_{\f O_0^A}}=\f O_0^A\ee 
by \eqref{seq}, this implies that    \be\l{yhb} \f O_2^X =\f O_0^A.\ee Finally, 
it follows from \eqref{yhb} by the second part of Lemma \ref{ser} that 
\be \l{rq3} \f O_2^{ \theta_{\f O_0^B}}\preceq \f O_1^{ \theta_{\f O_0^A}},\ee
implying that   $\f O_0^B$ is  non-ramified by \eqref{seq}.
 \qed

\bt \l{3} Let $A_1$, $A_2$ be non-special rational functions of degree at least two.
Then
 there exist rational functions $X_1$, $X_2$, $B_1$, $B_2$ such that  $X_1$, $X_2$ are Galois coverings of $\P^1(\C)$ by $\P^1(\C)$, $B_1$, $B_2$  are not generalized Latt\`es maps, 
the diagram 
\be  \l{ept} 
\begin{CD}  
 (\P^1(\C))^2 @>(B_1,B_2)>> (\P^1(\C))^2
\\ 
@V (X_1,X_2)  VV @VV  (X_1,X_2) V\\ 
 (\P^1(\C))^2 @>(A_1,A_2)>> (\P^1(\C))^2 
\end{CD}
\ee
commutes, and every irreducible $(A_1,A_1)$-periodic (resp. preperiodic) curve is the \linebreak
$(X_1,X_2)$-image of some irreducible $(B_1,B_2)$-periodic (resp. preperiodic) curve. 
\et
\pr Applying Theorem \ref{33} to $A_1$ and $A_2$, we obtain the commutative diagram
\be 
\begin{CD} 
(\P^1(\C))^2 @> (B_1,B_2)>>(\P^1(\C))^2 \\ 
@V (\theta_{\f O_0^{A_1}},\theta_{\f O_0^{A_2}}) VV @VV (\theta_{\f O_0^{A_1}},\theta_{\f O_0^{A_2}})  V\\ 
 (\P^1(\C))^2 @> (A_1,A_2)>> (\P^1(\C))^2, 
\end{CD} 
\ee
where $B_1,$ $B_2$ are not generalized Latt\`es map, and the use of  Lemma \ref{lk} finishes the proof. \qed

\br \l{rem2}
Note that in fact we proved a more precise version of Theorem \ref{3} with the concrete representation  
$$X_1=\theta_{\f O_0^{A_1}}, \ \ \ \ X_2=\theta_{\f O_0^{A_2}}$$  suitable for applications. 
\er

\end{subsection}

\end{section}

\begin{section}{Points with Zariski dense  orbits}
\begin{subsection}{Canonical heights and semiconjugacy} 

Let $K$ be a field of characteristic zero, which is assumed to be a subfield of $\C$ finitely generated over $\Q$, and   
$A\in K(z)$ a non-special rational function of degree $n\geq 2$.   
In this section, we give a criterion for the  $(A,A)$-orbit of a point $(x_0,y_0)\in (\P^1(\overline{K}))^2$ to be Zariski dense in $(\P^1(\C))^2$ in terms of canonical heights of $x_0$ and $y_0$. As an application, we prove a version of a conjecture of Zhang on the existence of rational points with Zariski dense forward orbits for endomorphisms 
of $(\P^1(K))^2$. We first assume that $K$ is a number field and use the Weil height. Then we explain how to extend our results to the general case using the Moriwaki height.    

Let $K$ be a number field. For $x\in \P^1(\overline{K})$ we denote by $h(x)$ the (logarithmic) Weil height of $x$.
We recall that  for any rational function  $R\in   \overline{K}(z)$ of degree $m$ there exists 
a constant $C_1>0$ depending only on $R$ such that for every $x\in \P^1(  \overline{K})$ the 
inequality 
\be \l{we} \vert h(R(x))-mh(x)\vert < C_1\ee holds. 
Furthermore, 
by the Northcott theorem, for any numbers $D_1,D_2 > 0$ there are only finitely many points $x\in \P^1(\overline{K})$ satisfying the conditions $$h(x) \leq D_1, \ \ \ \ \ [\Q(x):\Q]\leq D_2$$
(see e.g. \cite{silver}). 
 
Following \cite{cs}, for $A\in \overline{K}(z)$ of degree $n\geq 2$ we define the canonical height (associated to $A$) of 
a point   $x\in \P^1(\overline{K})$ as the limit 
\be \l{deff} \hat h_A(x)=\lim_{r\to \infty} \frac{h(A^{\circ r}(x))}{n^r}.\ee
We recall the following properties of the canonical height (see \cite{cs}, \cite{silver}). 
First, for every $x \in \P^1(\overline{K})$ the equality 
\be \l{kis} \hat h_A(A(x))=n\hat h_A(x)\ee holds. Second, there is a constant $C_2>0$ depending only on $A$ such that 
\be \l{to} |\hat h_A(x)-h(x)| < C_2\ee for
every  $x \in \P^1(\overline{K})$. Third, a point $x \in \P^1(\overline{K})$ is $A$-preperiodic if and only if $\hat h_A(x)=0.$ Finally, we mention that the function $\hat h_A: \P^1(\overline{K})\rightarrow \R$ is defined by the conditions \eqref{kis} and \eqref{to} in a unique way. 

Note that \eqref{we}  and  \eqref{to} imply that  
for any $R\in \overline{K}(z)$ of degree $m$ there exists a 
 constant $C_3>0$ depending only on $A$ and $R$ such that 
\be \l{zv} \vert \hat h_A(R(x))-m\hat h_A(x)\vert < C_3\ee 
 for every $x\in \P^1(\overline{K})$. Specifically,  
$$\vert \hat h_A(R(x))-m\hat h_A(x)\vert < 
\vert \hat h_A(R(x))-h(R(x))\vert +\vert h(R(x))-
m\hat h_A(x)\vert <$$
$$\vert \hat h_A(R(x))-h(R(x))\vert +\vert h(R(x))-mh(x)\vert +
\vert mh(x)-m\hat h_A(x)\vert <C_2+C_1+mC_2.$$

\bp \l{sehe} Let $A\in \overline{K}(z)$ and $B\in \overline{K}(z)$ be rational functions of degree at least two, and $X\in \overline{K}(z)$ a rational function of degree at least one such that the equality $A\circ X=X\circ B
$ holds.  Then for every point $x\in \P^1(\overline{K})$ the equality 
$$\hat h_A(X(x))=\hat h_B(x)\deg X$$ holds. 
\ep 
\pr Setting $n=\deg A=\deg B$ and using  inequality \eqref{we}, we have:   
$$\hat h_A(X(x))=\lim_{r\to \infty} \frac{h(A^{\circ r}(X(x)))}{n^r}=\lim_{r\to \infty} \frac{h(X(B^{\circ r}(x)))}{n^r}=$$
$$\lim_{r\to \infty} \frac{h(B^{\circ r}(x))\deg X+O(1)}{n^r}=\hat h_B(x)\deg X. \eqno{\Box}
$$

Notice that Proposition \ref{sehe} implies the following  known corollary.

\bc \l{kan} Let $A\in \overline{K}(z)$ be a rational function of degree $n\geq 2$ and $V\in \overline{K}(z)$ a rational function of degree $m\geq 2$ commuting with some iterate of $A$. Then 
 $\hat h_V=\hat h_A$. 
\ec 
\pr Proposition \ref{sehe} implies that if rational functions $B$ and $V$ commute, then 
$$\hat h_B(V(x))=\hat h_B(x)\deg V$$  for every $x\in \P^1(\overline{K}).$ Therefore, since the function $\hat h_A$ 
is defined by the conditions \eqref{kis} and \eqref{to} in a unique way,  for commuting $B$ and $V$ the equality $\hat h_V=\hat h_B$ holds. Thus, the condition of the corollary implies that 
$\hat h_V=\hat h_{A^{\circ l}}$ for some $l\geq 1.$ On the other hand, it follows from  \eqref{deff} that $\hat h_A=\hat h_{A^{\circ l}}$. \qed

\end{subsection} 
\begin{subsection}{\l{sec} Points with dense orbits: the case of a number  field}
Let $A\in \C(z)$ be a rational function of degree $n\geq 2$, and $n_0$ a minimum natural number such that $n=n_0^k$ for some $k\geq 1.$ Let us recall that by the Ritt theorem (\cite{r}) if a rational function $V\in \C(z)$  of degree $m\geq 2$ commutes  with $A$, then  $A$ 
and $V$ have a common iterate, unless they are both  special.  
Therefore, if $A$ is not special, there  
 exist $r,s\in \N$ such that 
\be \l{ka} V^{\circ r}=A^{\circ s},\ee  implying that  $m=n^{s/r}=n_0^l$ for some $l\in \N$.  
\vskip 0.2cm
\noindent{\it Proof of Theorem \ref{dh}.}  Let $(x_0,y_0)\in (\P^1(\overline{K}))^2$ be a point, and $\f O$ its $(A,A)$-orbit. Assume that 
the Zariski closure of $\f O$ in $(\P^1(\C))^2$ does not coincide with $(\P^1(\C))^2$.
It is easy to see (see e.g. \cite{ms}, Lemma 7.20) that then all but finitely many elements of $\f O$ are contained in some $(A,A)$-invariant   
algebraic set $Z\subset (\P^1(\C))^2$. Moreover, if  $x_0$ and $y_0$ are not  preperiodic points of $A$, then  $Z$ is a finite union of curves
 that are not vertical or horizontal lines.
Therefore, there exists an $(A,A)$-preperiodic curve $C\subset (\P^1(\C))^2$ that is not a vertical or horizontal line such that $(x_0,y_0)\in C.$  

Assume first that $A$ is not a generalized Latt\`es map. Then  Theorem \ref{2+} implies that $C$ is  
  a component of a separated variable curve $$V_1(x)-V_2(y)=0,$$ where  $V_1, V_2\in \C(z)$ are rational functions 
 commuting with some iterate of $A$. 
Moreover, we can assume that $V_1$ and $V_2$ belong to $\overline{K}(z)$ 
 (see Remark \ref{br3}). 
Hence, by Corollary \ref{kan}, \be \l{ol1} \hat h_{V_1}=\hat h_{V_2}= \hat h_A,\ee and, in addition,  \be \l{ol2} \deg V_1=n_0^{l_1}, \ \ \  \ \deg V_2=n_0^{l_2}\ee 
for some $l_1,l_2\in \N$. 

Since $V_1(x_0)=V_2(y_0)$, the equality 
$$\hat h_A(V_1(x_0))=\hat h_A(V_2(y_0))$$ holds.  
On the other hand, by \eqref{ol1} and \eqref{ol2}, we have:
$$\hat h_A(V_1(x_0))=\hat h_{V_1}(V_1(x_0))=n_0^{l_1}\hat h_{V_1}(x_0)=n_0^{l_1}\hat h_A(x_0),$$ 
$$\hat h_A(V_2(y_0))=\hat h_{V_2}(V_2(y_0))=n_0^{l_2}\hat h_{V_2}(y_0)=n_0^{l_2}\hat h_A(y_0).$$ 
Therefore,  equality \eqref{osn} holds for $l=l_1-l_2$.  

Assume now that  $A$ is a generalized Latt\`es map,  and 
let $(x_0,y_0)\in (\P^1(\overline{K}))^2$ be a point such  that 
 $x_0$ and $y_0$ are not preperiodic points of $A$, and the canonical heights of $x_0$ and $y_0$ do not satisfy condition \eqref{osn}.
By Theorem \ref{3}, 
 there exist rational functions $X$ and $B$ such that equality \eqref{i1} holds, $X:\P^1(\C)\rightarrow \P^1(\C)$ is a Galois covering, and  $B$  is not a generalized Latt\`es map. Moreover, in fact, $X=\theta_{\f O_0^A}$, implying that $c(\f O_2^X)=c(\f O_0^A)$ (see Remark \ref{rem2}). 
Let us observe that without loss of generality we can assume that $X$ and $B$ belong to $\overline{K}(x)$. Indeed, it  is well-known that for any Galois covering  $X:\P^1(\C)\rightarrow \P^1(\C)$ there exist M\"obius transformations $\delta_1,$ $\delta_2$ such that the function 
$\delta_1\circ X\circ \delta_2$  is ramified over $0,1,\infty$ and has rational coefficients. Since obviously  $c(\f O_2^A)\subset \overline{K}$ and $$ c(\f O_2^X)=c(\f O_0^A)\subseteq c(\f O_2^A)$$ by Lemma \ref{toch}
 (in case $\deg A<5$, we can consider $A^{\circ 3}$ instead of $A$), this implies 
that  for some M\"obius transformation $\delta$ the function $X'=X\circ \delta$ belongs to $\overline{K}(x)$. It follows now from \eqref{i1} that  preimages of any point $x\in \overline{ K}$ under the function $B'=\delta^{-1}\circ B\circ \delta$ belong to $\overline{ K}$, implying that $B'\in \overline{K}(x)$.

Let now $x_0'\in \P^1(\overline{ K})$ and $y_0'\in \P^1(\overline{ K})$ be arbitrary points such that the equalities $X(x_0')=x_0$ and $X(y_0')=y_0$ hold. It is easy to see that equality \eqref{i1} implies that 
 $x_0'$ nor $y_0'$ are not preperiodic points of $B$. Moreover, Proposition \ref{sehe}  implies that 
the canonical heights of $x_0'$ and $y_0'$ associated to $B$ do not satisfy the condition 
$$\hat h_B(y_0')=n_0^{l}\hat h_B(x_0'),\ \ \ \ \ \ \ l\in \Z.$$ 
Applying the already proved part of the theorem, we conclude that the $(B,B)$-orbit of $(x_0',y_0')$ is dense in  $(\P^1(\C))^2$. Since for any $(A,A)$-preperiodic curve  $\c C$ in $(\P^1(\C))^2$ the preimage 
 $\c E=(X,X)^{-1}(\c C)$ is a union of $(B,B)$-preperiodic curves, this implies that   $(A,A)$-orbit of $(x_0,y_0)$ is dense in  $(\P^1(\C))^2$. \qed 

\vskip 0.2cm
Let us recall that the Zhang conjecture about orbits states that if $\phi: X \rightarrow  X$ is a polarizable dynamical system
over some number field $K$, then there exists a point $a \in X(\overline{ K})$ whose forward
$\phi$-orbit is Zariski dense (see [20]). More generally, it was conjectured in the paper \cite{ms} that if $X$ is an irreducible variety
over an  algebraically closed field of characteristic
zero $K$ and $f : X \rightarrow X$ is a dominant rational self map such that 
there do not exist a positive dimensional algebraic
variety $Y$ and a dominant rational map $g : X \rightarrow Y$ for which $g \circ f = g$, then
there exists $a\in X(K)$ with a Zariski dense forward orbit. 

For a detailed discussion of the above ``Zariski dense orbit conjecture'' and a description of a few special cases in which it is known to be true we refer the reader to the recent paper \cite{k}. In particular, the addendum to \cite{k} contains a proof of the Zariski dense orbit conjecture for endomorphisms of $(\P^1(K))^2$, 
which is based on the main result of \cite{k} about the existence of Zariski dense orbits for endomorphisms of projective surfaces. 
 Below, we give an alternative proof of the Zariski dense orbit conjecture for endomorphisms of $(\P^1(K))^2$, which is based on Theorem \ref{dh}. Notice that  for endomorphisms induced by special rational functions  the truth  of the Zariski dense orbit conjecture follows from known results (see \cite{k} and also \cite{ms} for the polynomial case). Thus, as before, we will consider  only the case of non-special functions.

\bt \l{z1} Let $K$ be a number field and $A_1,A_2\in K(z)$ non-special rational functions of degree at least two. Then
there is a point in $(\P^1( K))^2$ whose $(A_1,A_2)$-forward orbit is Zariski dense in $(\P^1( K))^2$.
\et 
\pr Let us show first that  if $A_1=A_2=A\in \overline{K}(z)$, then there exists a point in $(\P^1(\Q))^2$ whose $(A,A)$-orbit is dense in $(\P^1(\C))^2$.

Let $R\in  \Q(x)$ be an arbitrary rational function  of degree two, say, $z^2$. By the Northcott theorem, there is a point  $x_0\in \P^1( \Q)$ with $\hat h_A(x_0)>C_3,$ 
where $C_3$ is a constant such that \eqref{zv} holds. 
Setting $n=\deg A$ and assuming that $n_0>2,$ we obtain from \eqref{zv} the inequalities  
$$\hat h_A(R(x_0))>2\hat h_A(x_0)-C_3>\hat h_A(x_0)$$ and 
$$\hat h_A(R(x_0))<2\hat h_A(x_0)+C_3<3\hat h_A(x_0)\leq n_0\hat h_A(x_0). $$
Thus, the point $y_0=R(x_0)\in \P^1( \Q)$ satisfies 
$$0<\hat h_A(x_0)<\hat h_A(y_0)<n_0\hat h_A(x_0),$$ 
implying that 
$x_0$ and $y_0$ are not preperiodic points of $A$, and the canonical heights of $x_0$ and $y_0$ do not satisfy condition \eqref{osn}. Therefore, the $(A,A)$-orbit of $(x_0,y_0)$ is dense in $(\P^1(\C))^2$ 
by Theorem \ref{dh}. 

In case $n_0=2$, instead of a rational function of degree two we can take any rational function $R\in  \Q(x)$ of degree three, obtaining 
 $$\hat h_A(R(x_0))>3\hat h_A(x_0)-C_3> 2\hat h_A(x_0)$$ and 
$$\hat h_A(R(x_0))<3\hat h_A(x_0)+C_3<4\hat h_A(x_0),$$
so that the point $y_0=R(x_0)\in \P^1( \Q)$ satisfies 
$$0<n_0\hat h_A(x_0)<\hat h_A(y_0)<n_0^2\hat h_A(x_0).$$ 

Assume now that  $A_1, A_2\in K(z)$ are arbitrary non-special rational functions of degree at least two, and let $(x_0,y_0)\in(\P^1( K))^2$ be an arbitrary point. If  the 
$(A_1,A_2)$-orbit of $(x_0,y_0)$ is dense in $(\P^1( K))^2$, we are done. Otherwise, infinitely many elements of the orbit belong to an $(A_1,A_2)$-periodic curve $C$ defined over $K$.  By Theorem \ref{in1}, there exist 
rational functions $X_1,X_2,B\in \overline{K}(x)$ and $r\geq 1$ such that the diagram 
$$
\begin{CD} 
(\P^1(\C))^2 @>(B,B)>>(\P^1(\C))^2 \\ 
@V (X_1,X_2)  VV @VV  (X_1,X_2) V\\ 
 (\P^1(\C))^2 @>(A_1^r,A_2^r)>> (\P^1(\C))^2
\end{CD}
$$
commutes, and $X_1,X_2$ parametrize $C.$ Moreover, since $C$ contains infinitely many points with coordinates in $\P^1( K)$, there exists a parametrization of $C$ defined over $K$. Thus, without loss of generality we may assume that  $X_1, X_2\in K(x)$. 
By what is proved above there is  a point $(x_0,y_0)\in(\P^1(\Q))^2$ whose  $(B,B)$ orbit  is dense in $(\P^1(\C))^2$,   implying  that the $(A_1^r,A_2^r)$-orbit of the point 
$(X_1(x_0),X_2(y_0))\in(\P^1( K))^2$ is dense in $(\P^1(\C))^2$ and hence in $(\P^1( K))^2$.  \qed 

\end{subsection} 
\begin{subsection}{Points with dense orbits: the case of an arbitrary  field}
 Let us recall that by the results of Moriwaki (\cite{mor}) for every  field $K$  finitely generated over $\Q$ one can define the height function $\mathfrak{h}$ on $\P^1(\overline{K})$ satisfying the following two properties used in Section \ref{sec}. For any  $R\in   \overline{K}(x)$ of degree $m$ there exists 
a constant $C_1>0$ such that for every $x\in \P^1(  \overline{K})$ the equality 
$$ \vert \mathfrak{h}(R(x))-m\mathfrak{h}(x)\vert < C_1$$ holds, and for any  $D_1,D_2 > 0$ there are only finitely many points $x\in \P^1(\overline{K})$ satisfying the conditions $$\mathfrak{h}(x) \leq D_1, \ \ \ \ \ [\Q(x):\Q]\leq D_2.$$ Defining now the canonical height corresponding to  a rational function 
$A\in \overline{K}(x)$ by the formula 
$$\hat{\mathfrak{h}}_A(x)=\lim_{r\to \infty} \frac{\mathfrak{h}(A^{\circ r}(x))}{n^r}$$ 
and repeating verbatim the proofs of Proposition \ref{sehe} and Theorem \ref{dh}, we obtain the following result.

\bt \l{dh++} Let $K\subset \C$ be a field finitely generated over $\Q$ and  $A$ a non-special rational function of degree
 $n\geq 2$ defined over $K$. 
Then the $(A,A)$-orbit of a point $(x_0,y_0)\in (\P^1(\overline{K}))^2$ is Zariski dense in $(\P^1(\C))^2$, unless 
either $x_0$ or $y_0$ is a preperiodic point of $A$, or the canonical heights of $x_0$ and $y_0$ satisfy the condition 
$$ \hat{\mathfrak{h}}_A(y_0)=n_0^{l}\hat{\mathfrak{h}}_A(x_0),\ \ \ \ \ \ \ l\in \Z,$$ where $n_0$ is a minimum natural number such that $n=n_0^k$ for some $k\geq 1.$ \qed 
\et 

In turn, Theorem \ref{dh++} implies the following statement.

\bt \l{z2} Let $K\subset \C$ be a  field  and $A_1,A_2\in K(z)$ non-special rational functions of degree at least two. Then
there is a point in $(\P^1( K))^2$ whose $(A_1,A_2)$-forward orbit is Zariski dense in $(\P^1( K))^2$. 
\et 
\pr Defining $K'$ as the subfield of $K$ generated by the coefficients of $A_1$, $A_2$ and 
arguing as in the proof of Theorem \ref{z1} we can find a point in  $(\P^1( K'))^2$ whose 
$(A_1,A_2)$-orbit is dense in $(\P^1(\C))^2$ and hence in $(\P^1( K))^2$. \qed

\end{subsection}

\end{section}

\begin{section}{Finiteness theorems}
\subsection{Formulation of results}
In this section, we prove several results, which can be considered as quantitative analogues of 
 results of the paper \cite{aol} in a slightly simplified setting. As an application, we prove Theorem \ref{5} and the finiteness 
of the number of $(A_1,A_2)$-invariant curves of any given bi-degree $(d_1,d_2)$.

The first result is following.

\bt \l{6} 
There exists a function $\phi:\mathbb N\times \mathbb N\rightarrow \R$ with the following property. For
any non-special  rational function $A$ of degree at least two and rational function $X$ such that  for every $d\geq 1$
the algebraic curve \be \l{cura} A^{\circ d}(x)-X(y)=0\ee has a factor of genus zero, there exists  
$N\leq\phi(\deg A,\deg X)$ such that the equality \be \l{puzo0} A^{\circ N}\circ \theta_{\f O_0^A}=X\circ R\ee holds for some rational function $R.$  
\et

Note that the assumption of Theorem \ref{6} about curves \eqref{cura} holds for any pair of rational functions $A$ and $X$ satisfying \eqref{ssee} for some rational function $B$.  Indeed, it follows from \eqref{ssee} that 
$$A^{\circ d}\circ X=X\circ B^{\circ d}, \ \ \ \ \  \ d\geq 1,$$ implying that  the curve \eqref{cura} has a component of genus zero with the parametrization $t\rightarrow (X(t),B^{\circ d}(t))$. 
Similarly, 	the above assumption  holds for any $A$ and $X$ satisfying \eqref{decc}. However, in this case 
Theorem \ref{5} provides a more precise conclusion which permits to get rid of the function $\theta_{\f O_0^A}$ in \eqref{puzo0}. On the other hand, if $A$ is not a generalized Latt\`es map, then 
$\theta_{\f O_0^A}$ reduces to the identical map even in the more general setting of Theorem \ref{6}.

 We say that two rational functions  $W_1$ and $W_2$ are $\mu$-{\it equivalent} if there exists a M\"obius transformation $\mu$ such that $$ W_1=W_2\circ \mu.$$
The next result is a weaker form of Theorem \ref{6}, which holds, however, for all functions $A$ including special, 
for which the function $\theta_{\f O_0^A}$ is transcendental or is not defined.

\bt \l{7} 
There exists a function $\chi:\mathbb N\times \mathbb N\rightarrow \mathbb R$ with the following property. For any rational functions $A$ of degree $m\geq 2$ and integer $n\geq 1$, 
there exist at most $\chi(m,n)$ classes of $\mu$-equivalence of rational functions $X$ of degree $n$ 
 such that 
for every $d\geq 1$
the algebraic curve $$ A^{\circ d}(x)-X(y)=0$$ has a factor of genus zero.
\et

Let $A$ be a rational function. 
We denote by $D= D\Big[A,N,\langle W_d,h_d\rangle_{d=1}^N\Big]$
a  commutative diagram  of the form 
\be \l{coma}
\begin{CD}
\P^1(\C)  @> h_N >>  \P^1(\C)  @. \hskip 0.6cm \dots \hskip 0.6cm \P^1(\C) @> h_2 >>  \P^1(\C) @>h_1 >> \P^1(\C) \\
@VV W_N V @VV W_{N-1} V  \hskip 1.2cm @VV  W_2 V  @VV W_1 V @VV W_0 V \\
\P^1(\C)  @> A >>   \P^1(\C) @. \hskip 0.6cm \dots \hskip 0.6cm \P^1(\C)\ @> A >> \P^1(\C) @> A >> \P^1(\C),
\end{CD}
\ee 
where $h_d,W_d,$ $1\leq d \leq N,$ and $W_0$ are  rational functions. 
We say that $D$ is {\it good} if 
 for any $d_1,d_2$,  $0\leq d_1<d_2\leq N$,  the functions
\be \l{mapa} W_{d_1}, \ \ \ h_{d_1+1}\circ h_{d_1+2}\circ \dots \circ h_{d_2}, \ \ \ A^{\circ(d_2-d_1)}, \ \ \ W_{d_2}\ee form a good solution of equation \eqref{m}. Note that if $D$ is good, then by Lemma \ref{good} $\deg W_d=\deg W_0,$ $ d\geq 1.$   
For a good diagram $D$, we set $$m_D=\deg A, \ \ \ \ n_D=\deg W_0.$$ We call the number $N$ the {\it length} of $D$. 
For a diagram $D=D\Big[A,N,\langle W_d,h_d\rangle_{d=1}^N\Big]$ and $j_1,$ $j_2,$ $0\leq j_1<j_2\leq N$,  we denote by $D_{j_1,j_2}$ the sub-diagram of $D$ bounded by the arrows $W_{j_1}$ and   $W_{j_2}$.

Let $r$, $1\leq r\leq N$, be an integer. We say that $D=D\Big[A,N,\langle W_d,h_d\rangle_{d=1}^N\Big]$ is $r$-{\it periodic}   if 
for every $j,$ $0\leq j\leq N-r$, the equality 
$$ W_{j+r}=W_{j}\circ \alpha_{j}$$ holds for some M\"obius transformation $\alpha_{j}$.
We say that $D$ is {\it periodic} if it is $r$-periodic for some $r,$ $1\leq r\leq N$.
Finally, we say that $D$ is {\it preperiodic} if for some $N_0,$ $0\leq N_0\leq N-1,$ the sub-diagram $D_{N_0,N}$ is periodic. 

The last of the analogues of results of the paper \cite{aol} proved in this section  
is following.

\bt \l{xyi} There exists a function $\psi:\mathbb N \times\mathbb N \rightarrow \R$ with the following property. Any  good diagram $D=D\Big[A,N,\langle W_d,h_d\rangle_{d=1}^N\Big]$ such that $m_D\geq 2$  
and \linebreak  $N> \psi(m_D,n_D)$  is preperiodic.
\et

\subsection{Proof of Theorem \ref{xyi}}

As in \cite{aol}, we use the following result proved in \cite{cur}.

\bt \l{m2} 
Let $U$ be a rational function of degree $n$.
Then for any rational function 
$V$ of degree $m$ such that the curve $ \f E_{U,V}:\,
U(x)-V(y)=0$ is irreducible
the inequality 
\be \l{ma} 
g(\f E_{U,V}) >\frac{m-84n+168}{84}
\ee holds, unless $\chi(\f O_2^U)\geq 0$. \qed
\et 

We also  need the following lemma, which is a particular case of Theorem 2.4 in the paper \cite{aol}.

\bl \l{xriak+} Let $R$ be a compact Riemann surface, $f:\, R\rightarrow \P^1(\C)$ a holomorphic map, and  $\f O$ an orbifold. 
Then \be \l{dfc} \theta_{\f O}=f\circ h\ee for some holomorphic map $h:\t{\f O}\rightarrow R$  if and only if $\f O^f_2\preceq \f O$. \qed
\el
For brevity, 
we call a rational function $f$ satisfying \eqref{dfc} a {\it compositional left factor} of  $\theta_{\f O}$. More precisely, by  
a compositional left factor of a holomorphic map \linebreak $f:\, R_1\rightarrow R_2$ between Riemann surfaces, 
we mean any holomorphic map \linebreak  $g:\, R'\rightarrow R_2$ between Riemann surfaces such that $f=g\circ h$ for some  holomorphic map  $h:\, R_1\rightarrow R'.$

\bl \l{lu} There exists a function $\kappa:\, \N\rightarrow \N$ with the following property. For any orbifold $\f O$ with $\chi(\f O)\geq 0$ 
there exist at most $\kappa(n)$ classes  of $\mu$-equivalence  of rational functions  $f$ of degree $n$  with $\f O_2^f=\f O$.
\el 
\pr 
By Lemma \ref{xriak+}, the equality $\f O_2^f=\f O$ implies that $f$ is a compositional left factor of $\theta_{\f O}$.  Moreover, it is easy to see that the equality 
\be \l{krot} \theta_{\f O_2^{f}}=f\circ \theta_{\f O_1^{f}}\ee holds. 
Therefore, for any fixed $\f O$, the number of $\mu$-equivalence classes of rational functions  $f$ of degree $n$  with $\f O_2^f=\f O$ does not exceed the number of subgroups of index $n$ in the group $\Gamma_{\f O}$. On the other hand, the number of such subgroups  can be bounded in terms of $n$, since  $\Gamma_{\f O}$  is finitely generated. Since the group $\Gamma_{\f O}$ is defined by the signature of $\f O$, this implies that to prove the lemma we only must show that for orbifolds whose signatures belong to the 
{\it infinite} series  $\{l,l\},$ $l\geq 2,$ and $\{2,2,l\},$ $l\geq 2,$ from the lists \eqref{list}, \eqref{list2}, the bounds for the number of $\mu$-equivalence  classes are uniform.

It is well-known (see e.g. 
Corollary 2.7 in \cite{aol}) that if 
$\f O_2^f$ is defined by the conditions
\be \l{los1} \nu_2^f(0)=l, \ \ \ \ \nu_2^f(\infty)=l,\ee then 
\be \l{koz1} f=z^{l}\circ \mu\ee 
for some M\"obius transformation $\mu$, while if $\f O_2^f$ is defined by the conditions
\be \l{los2} \nu_2^f(-1)=2, \ \ \ \ \nu_2^f(1)=2, \ \ \ \  \nu_2^f(\infty)=l, \ee  
then either \be \l{koz2} f=\frac{1}{2}\left(z^l+\frac{1}{z^l}\right)\circ \mu,\ee
or \be \l{koz3} f=\pm T_l\circ \mu\ee for some M\"obius transformation $\mu.$ Therefore, 
there exists exactly one $\mu$-equivalence class of
rational functions $f$ of degree $n$ such that the signature of $\f O_2^f$ 
belongs to the series  $\{l,l\},$ $l\geq 2,$ and there exist at most three such classes if  the signature of 
$\f O_2^f$ 
belongs to the series $\{2,2,l\},$ $l\geq 2.$ \qed

\bl \l{nova} Let $D=D\Big[A,N,\langle W_d,h_d\rangle_{d=1}^N\Big]$ be a diagram 
such  that 
\be \l{ragna} \C(h_d,W_d)=\C(z), \ \ \ 1\leq d \leq N.\ee Assume that 
\be \l{usl} W_{r}=W_{0}\circ \mu\ee for some $r,$ $1\leq r\leq N,$ and  M\"obius transformation $\mu$.
Then $D$ is good and $r$-periodic.
\el 
\pr Since \eqref{ragna} implies that 
the map $$t\rightarrow (h_d(t),W_d(t)),\ \ \ 1\leq d\leq N,$$ is a generically one-to-one parametrization of some 
component of the curve $$W_{d-1}(x)-A(y)=0,$$ we  see that 
\be \l{asas} \deg W_N\leq \deg W_{N-1}\leq \dots \leq \deg W_1\leq \deg W_0.\ee
Thus, \eqref{usl} yields that 
$$\deg W_{r}=\deg W_{r-1}=\dots = \deg W_{1}= \deg W_{0},$$ implying 
by Lemma \ref{good}  and Theorem \ref{sum+} that the sub-diagram $D_{0,r}$   is  good. In particular, the fiber product of $W_{0}$ and $A$ has a unique component and the functions  $W_1,$ $h_1$ are defined 
by $W_{0}$ in a unique way up to natural isomorphisms. It follows now from \eqref{usl} that  
the fiber product of $W_{r}$ and $A$ also has a unique component and  $$W_{r+1}=W_{1}\circ \mu'$$ for some M\"obius transformation $\mu'$. In particular, the sub-diagram $D_{0,r+1}$   is  good.
 Continuing arguing in this way, we conclude that $D$ is good and $r$-periodic. \qed

\vskip 0.2cm

\noindent{\it Proof of Theorem \ref{xyi}.}
We first prove the theorem under 
the additional assumption
\be \l{bur} \chi(\f O_2^{W_d})\geq 0,\ \ \ \ 0\leq d\leq N.\ee 
 For a good diagram $D=D\Big[A,N,\langle W_d,h_d\rangle_{d=1}^N\Big]$ define $k=k(D)$ as the number of {\it distinct} orbifolds among the orbifolds $\f O_2^{W_{d}},$ $0\leq d\leq N.$ To prove the theorem it is enough to show that there exists a function $C=C(m_D)$ such that 
\be \l{inde} k(D)\leq  C(m_D).\ee Indeed, if \eqref{inde} holds, then 
 Lemma \ref{lu} and the box principle imply that 
whenever $$N> \psi=C(m_D)\kappa(n_D),$$ 
there exist $j_1,$ $j_2,$ $0\leq j_1<j_2\leq N$, such that  $W_{j_2}$ and $W_{j_1}$ are $\mu$-equivalent. 
Since equalities \eqref{ragna} hold by Lemma \ref{good},  this implies by Lemma \ref{nova} that the sub-diagram $D_{j_1,N}$ is $(j_2-j_1)$-periodic.

To prove \eqref{inde} it is enough to  bound  in terms of $m_D$ 
the number of distinct sets  among  the sets $c(\f O_2^{W_{d}}),$  $0\leq d\leq N,$ and 
the number of distinct signatures  among  the signatures
$\nu(\f O_2^{W_{d}}),$ $0\leq d\leq N.$ 
Since $$ A:\f O_2^{W_{d+1}}\rightarrow \f O_2^{W_{d}}, \ \ \ \ 0\leq d \leq N-1,$$  
is a minimal holomorphic map between orbifolds by Theorem \ref{goodt}, it follows from  Lemma \ref{toch} that if $m_D>4$, then 
every set $c(\f O_2^{W_{d}}),$ $0\leq d\leq N-1,$ is a subset of the set $c(\f O_2^A)$. Since a rational function of degree $m$ has at most $2m-2$ critical values, this
 implies  that  the  number of distinct sets among 
the sets  $c(\f O_2^{W_{d}}),$ $0\leq d\leq N,$ is bounded in terms of $m_D$.
Moreover, this is  also true if  $m_D\leq 4$. Indeed, the inequality $m_D\geq 2$ implies 
the inequality $m_D^{\circ 3}>4$, and hence 
every set $c(\f O_2^{W_{d}}),$ $0\leq d\leq N-3,$ is a subset of the set $c(A^{\circ 3})$, since 
$$A^{\circ 3}:\f O_2^{W_{d+3}}\rightarrow \f O_2^{W_{d}}, \ \ \ \ 0\leq d \leq N-3,$$ 
 also are  minimal holomorphic maps.
Finally,  possible  signatures of the orbifolds $\f O_2^{W_{d}},$ $0\leq d\leq N,$
are contained in the lists \eqref{list}, \eqref{list2}, and by formulas \eqref{koz1}, \eqref{koz2}, \eqref{koz3}, if  $\nu(\f O_2^{W_{d}})=\{l,l\}$, $l\geq 2,$ then $l=n_D$, while if 
 $\nu(\f O_2^{W_{d}})=\{2,2,l\}$, $l\geq 2,$ then either
$l=n_D$ or $l=n_D/2$. Thus, 
 the   number of  distinct signatures  among  the signatures $\nu(\f O_2^{W_{d}}),$ $d\geq 0,$ does not exceed ten.

The proof of the theorem in the general case reduces to the case where \eqref{bur} is satisfied. 
Indeed, since the commutativity of diagram \eqref{coma} implies that the curves 
$$A^{\circ d}(x)-W_0(y)=0, \ \ \ 1\leq d\leq N,$$
have genus zero, applying Theorem \ref{m2} for $U=W_0$ and $V=A^{\circ N}$,  we see 
 that  whenever $$m_D^N>84(n_D-2)$$ the inequality $ \chi(\f O_2^{W_0})\geq 0$ holds. More generally, 
setting $U=W_i$, $0\leq i \leq N_0,$ and $V=A^{\circ (N-i)}$,  we see that 
 whenever \be \l{kroi} m_D^{N-N_0}>84(n_D-2),\ \ \   N_0\geq 0,\ee 
the inequalities 
$$ \chi(\f O_2^{W_d})\geq 0,\ \ \ \ 0\leq d\leq N_0,$$
 hold. 
Therefore, if  \be \l{funn} N>\psi=\log_{m_D}(84(n_D-2))+C(m_D)\kappa(n_D)+1 ,\ee then
the inequalities 
$$ \chi(\f O_2^{W_d})\geq 0,\ \ \ \ 0\leq d\leq C(m_D)\kappa(n_D)+1 ,$$
 hold. By the already proved part of the theorem, we conclude that  
there exist $j_1,$ $j_2,$ $0\leq j_1<j_2\leq C(m_D)\kappa(n_D)+1$, such that  $W_{j_2}$ and $W_{j_1}$ are $\mu$-equivalent, implying as above that $D$ is preperiodic. 
\qed

\vskip 0.2cm

\subsection{Proof of Theorem \ref{6}, Theorem \ref{7}, and Theorem \ref{5}}

\vskip 0.2cm
\noindent{\it Proof of Theorem \ref{6}.} Since for any holomorphic map $f:R\rightarrow R'$ between compact 
Riemann surfaces the inequality $g(R)\geq g(R')$ holds,  it follows from the universality property of the fiber product that  if for every $d\geq 1$
curve \eqref{cura}  has a factor of genus zero, then 
for every $N\geq 1$ there exists a  diagram $D$ of the form \eqref{coma} such that $W_0=X$ and the conditions   \eqref{ragna}, \eqref{asas} hold. 

Assume that for some $ l_1,$ $ l_2,$ $0\leq  l_1< l_2 \leq N,$ the condition 
\be \l{zaiaz} \deg W_{ l_2}=\dots = \deg W_{ l_1+1}= \deg W_{ l_1}\geq 2\ee holds. 
Then we conclude as in Lemma \ref{nova}  that 
the sub-diagram $D_{l_1,l_2}$  
is  good, and  applying  
Theorem \ref{xyi} to the diagram $D_{l_1,l_2}$ we see that either 
there exist $j_1,$ $j_2,$  $ l_1\leq j_1<j_2\leq  l_2$ such that 
$W_{j_2}$ and $W_{j_1}$ are $\mu$-equivalent,
or 
$$  l_1- l_2\leq \psi(\deg A,\deg W_{ l_1}),$$ implying that 
\be \l{bou}  l_1- l_2\leq \psi(\deg A,\deg X),\ee 
since the function in the right part of \eqref{funn} is increasing in the argument $n_D$.
It follows now from \eqref{asas} and \eqref{bou} 
that whenever \be \l{sati} N> \phi(\deg A,\deg X)=\psi(\deg A,\deg X)\cdot(\deg X-1)+1,\ee
 either \be \l{qew} \deg W_{N}=1,\ee 
or there exist a M\"obius transformation $\mu$ and integers $j_1,$ $j_2,$  $ 1\leq j_1<j_2\leq  N,$
such that  \be \l{ghf} W_{j_2}=W_{j_1}\circ \mu\ee 
and \be \l{kuku} \deg W_{j_2}=\deg W_{j_1}\geq 2.\ee

In the first case, the function
\be \l{by} R_1=h_1\circ h_2\circ \dots\circ h_{N}\circ W_{N}^{-1} \ee satisfies 
\be \l{pll} A^{\circ N}=X\circ R_1,\ee implying that 
$$A^{\circ N}\circ \theta_{\f O_0^A}=X\circ (R_1\circ \theta_{\f O_0^A}).$$
In the second case, the equality 
 \be \l{puzo1} A^{\circ j_1}\circ W_{j_1}=X\circ R_2\ee holds for
 the function  $$R_2 = h_1\circ h_2\circ \dots\circ h_{j_1}.$$
Furthermore,  since $D_{j_1,j_2}$  is good, it follows from \eqref{ghf} by Theorem \ref{goodt} that
$$A^{\circ (j_2-j_1)}:\f O_2^{W_{j_1}}\rightarrow \f O_2^{W_{j_1}}$$ is a minimal holomorphic map, and hence 
\be \l{70} \f O_2^{W_{j_1}}\preceq \f O_0^A,\ee by Theorem \ref{uni}. It follows now from Lemma \ref{xriak+} that the equality 
$$\theta_{\f O_0^A}=W_{j_1}\circ T$$ holds for some rational function $T$, implying by \eqref{puzo1} that 
$$A^{\circ j_1}\circ \theta_{\f O_0^A}=X\circ R_2\circ T.$$  
Thus, 
$$A^{\circ N}\circ \theta_{\f O_0^A}=X\circ R_2\circ T\circ F^{\circ (N-j_1)},$$ 
where $F$ is a rational function, which makes the diagram \eqref{dia3} commutative. \qed  
\vskip 0.2cm
\noindent{\it Proof of Theorem \ref{7}.}
We recall that if $R$ is a compact Riemann surface and $f:R\rightarrow \P^1(\C)$ is a holomorphic map, then 
functional decompositions  $f=U\circ V$, where $V:R\rightarrow R'$ and $U:R'\rightarrow \P^1(\C)$ are 
holomorphic maps between compact Riemann surfaces,  considered up to the equivalence 
\be \l{equi} U\rightarrow U\circ \mu, \ \ \ \ V\rightarrow\mu^{-1} \circ V,\ \ \ \ \mu\in \Aut(R'),\ee are in a 
one-to-one correspondence with imprimitivity systems of the monodromy group of $f$. Thus,  Theorem \ref{6} implies that for non-special $A$ the number of $\mu$-equivalence classes of rational functions $X$ of degree $n$ such that for every $d\geq 1$
the algebraic curve \eqref{cura} has a factor of genus zero
is bounded by the number of imprimitivity systems in the monodromy group of the function $A^{\circ N}\circ \theta_{\f O_0^A}.$ In turn, this number is bounded in terms of $m$ and $n.$ 

Assume now that $A$ is a Latt\`es map. In this case, it is still true that if $N$ satisfies \eqref{sati}, 
then either conditions \eqref{qew}  and \eqref{pll}, or  conditions \eqref{puzo1}  and \eqref{70} hold. 
Moreover,  $$\deg W_{j_1}\leq \deg W_0=n,$$ and \eqref{70} implies that $$\chi(\f O_2^{W_{j_1}})\geq \chi(\f  O_0^{A}) = 0.$$ 
It follows now from Lemma \ref{lu} that 
 the considered number of $\mu$-equivalence classes  is bounded by the total number of imprimitivity systems in the monodromy groups of a {\it finite number} of rational functions of the form $A^{\circ j_1}\circ W,$ where  $\deg W\leq n$, 
 $\f O_2^{W}\preceq \f O_0^A$, and     $j_1 \leq N.$

Finally, by Theorem 3.6 of \cite{aol}, if $A$ is conjugate to $z^{\pm m}$, then any $X$ satisfying the conditions of the theorem has the form   
$X=z^{n}\circ \mu$ for some $\mu\in \Aut(\P^1(\C))$, while  if $A$  is conjugate to 
 $\pm T_m$, then either $X=\pm T_{n}\circ \mu,$ or $$X=\frac{1}{2}\left(z^{n/2}+\frac{1}{z^{n/2}}\right)\circ \mu,$$ for some $\mu\in Aut(\P^1(\C))$. Thus, the theorem is true also in this case.
\qed

\vskip 0.2cm
\noindent{\it Proof of Theorem \ref{5}.} 
It follows from equality \eqref{decc}  
that the map \be \l{pa}t\rightarrow (A^{\circ (d-1)}(t),R(t))\ee is a parametrization of some irreducible
component of the curve $$A(x)-X(y)=0.$$ This parametrization is not necessary one-to-one. However, we can find a parametri\-zation  $W_1,h_1$ such that $\C(W_1,h_1)=\C(z)$. Moreover,  the functions $W_1,h_1$ satisfy  the equalities 
\be \l{e2}  A^{\circ (d-1)}(t)=W_1\circ H_1,  \ \ \ R=h_1\circ H_1\ee
for some rational function $H_1$.
In particular, the diagram 
$$
\begin{CD}
\P^1(\C) @> H_1 >>  \P^1(\C) @>h_1 >> \P^1(\C) \\
 @VV  z V  @VV W_1 V @VV X V \\
\P^1(\C)\ @> A^{\circ (d-1)} >> \P^1(\C) @> A >> \P^1(\C)
\end{CD}
$$ 
commutes.
Similarly,  the map $$t\rightarrow (A^{\circ (d-2)}(t),H_1(t))$$ is a parametrization of some irreducible
component of the curve $$A(x)-W_1(y)=0,$$ implying  that  
there exist rational functions $W_2,h_2$ and $H_2$ such that the equalities $$ A^{\circ (d-2)}(t)=W_2\circ H_2,  \ \ \ H_1=h_2\circ H_2, \ \ \ \C(W_2,h_2)=\C(z)$$ hold  
and the diagram 
\be 
\begin{CD}
\P^1(\C)  @> H_2 >>  \P^1(\C)  @> h_2 >>  \P^1(\C) @>h_1 >> \P^1(\C) \\
@VV z V @VV W_2 V    @VV W_1 V @VV X V \\
\P^1(\C)  @> A^{\circ (d-2)} >>    \P^1(\C)\ @> A >> \P^1(\C) @> A >> \P^1(\C),
\end{CD}
\ee commutes.
Continuing  arguing in the same way, for every $N\leq d$  we obtain diagram \eqref{coma}, such that $$W_0=X, \ \ \ \ W_N=z,$$
and the conditions   \eqref{ragna}, \eqref{asas} hold. 

Now, as in the proof of Theorem \ref{6} and Theorem \ref{7}, we conclude that if $N$ satisfies \eqref{sati}, then either equalities \eqref{qew}, \eqref{pll} hold, or there exist integers $j_1,$ $j_2,$  $ 1\leq j_1<j_2\leq  N,$
such that \eqref{ghf} and \eqref{kuku} hold. 
However, the last case is impossible. Indeed, if \eqref{ghf} holds, then 
Lemma \ref{nova} applied to the diagram $D_{j_1,d}$ 
implies that $$\deg W_d=\deg W_{ j_1},$$  
in contradiction with the conditions  $$\deg W_d=1,\ \ \ \ \deg W_{ j_1}\geq 2.\eqno{\Box}$$

\subsection{Finiteness 
of the number of invariant curves of a given bi-degree}

Let 
$R$ be a compact Riemann surface of genus zero or one, and $B:R\rightarrow R$ a   holomorphic map.
We denote by $G_1(B)$ the subgroup of $ \Aut(R)$  consisting of $\mu\in \Aut(R)$ such that $$B\circ \mu =B,$$ and   by $G_2(B)$  the subgroup consisting of $\mu$ such that $$\mu^{-1}\circ B\circ \mu =B.$$ 

\bl \l{kroo} The group  $G_1(B)$ is finite, and its order can be bounded in terms of the degree of $B$. The same conclusion holds for the group $G_2(B)$ whenever the degree of $B$ is at  least two.  
\el 
\pr  Assume first that $g(R)=0$, so that $B$ is a rational function and elements of $G_1(B)$ and $G_2(B)$ are M\"obius transformations. If $\deg B=1$, then the group $G_1(B)$ is trivial. So,  assume that $\deg B\geq 2.$ Let us observe that  
any $\mu\in G_1(B)$ permutes preimages of $(B^{\circ k})^{-1}(z_0)$ for any $z_0\in \P^1(\C)$ and $k\geq 1.$ Since each M\"obius transformation is determined by specifying its value at three distinct points,  this implies that the group $G_1(B)$ is finite and its order can be bounded in terms of $\deg B.$ 
 Similarly, 
any $\mu\in G_2(B)$ permutes $B$-periodic points of any given period $k\geq 1$,  
implying that the group $G_2(B)$ is finite.

If $g(R)=1$, then  any $\mu\in G_1(B)$ still permutes preimages of $(B^{\circ k})^{-1}(z_0)$, while any  $\mu\in G_2(B)$ permutes $B$-periodic points. 
Furthermore, any $\mu\in \Aut(R)$ is induced by a linear map 
$$F=\omega z+c,\ \ \ \ \omega, c\in \C,$$ where $\omega$ is an $l$th root of unity with $l=1,2,3,4,$ or $6$. Such  $\mu$ has $|\omega-1|^2$ fixed points, implying that  it is determined by its
values at $|\omega-1|^2+1$ distinct points. Thus, the same argument as above shows the finiteness of $ G_1(B)$
 and $G_2(B)$. \qed

\bl \l{bug} 
Let $A$ be a rational function of degree at least two, $R$  a compact Riemann surface of genus zero or one, and $X:R\rightarrow \P^1(\C)$ a holomorphic map. Then the number of   holomorphic maps $B:R\rightarrow R$ such that the diagram 
\be \l{krol} 
\begin{CD}
R @> B >> R\\
@VV X V @VV X V\\ 
\P^1(\C) @>A >>\P^1(\C)
\end{CD}
\ee
commutes is finite and can be bounded in terms of degrees of $A$ and $X$.
\el 
\pr Setting $F=A\circ X$, we see that any two functions $B$ and $B'$ making diagram \eqref{krol} commutative satisfy the equality 
$$F=X\circ B=X\circ B'.$$ Since the number of imprimitivity  systems in
 the monodromy group of $F$ is finite, this implies that 
there exist holomorphic maps $B_1,B_2,\dots, B_N:R\rightarrow R$  such that the equality $F=X\circ B$ 
 holds  for a  holomorphic map $B:R\rightarrow R$ if and only if there exists $\mu\in \Aut(R)$ such that 
\be \l{poil} X=X\circ \mu, \ \ \ \ \ B=\mu^{-1} \circ B_j\ee for some $j,$ $1\leq j \leq N.$  Moreover, the number $N$ is bounded in terms of degrees of $A$ and $X$, since
$\deg F=\deg A\cdot\deg X$. Finally, the number of  $\mu$ satisfying the first equality in \eqref{poil} is also bounded by  Lemma \ref{kroo}.  
\qed

\bt \l{4}  Let $A_1$, $A_2$ be  rational functions of degree $m\geq 2.$
Then for any pair of positive integers $(d_1,d_2)$ there exist at most finitely many $(A_1, A_2)$-invariant curves of bi-degree $(d_1,d_2).$ 
Moreover, there exists a function $\gamma:\mathbb N\times \mathbb N\times \mathbb N\rightarrow \mathbb R$
such that the number of these curves does not exceed $\gamma(m,d_1,d_2).$ 
\et
\pr
 Assume first that $A_1$, $A_2$ are not both Latt\`es maps. Then  
by Theorem \ref{in1} any irreducible invariant curve $\f C$  of bi-degree $(d_2,d_1)$ has genus zero and can be parametrized by rational functions 
$X_1$ and $X_2$ of degrees $d_1$ and $d_2$  correspondingly making the diagram 
\be \l{doi} 
\begin{CD} 
(\P^1(\C))^2 @>(B,B)>>(\P^1(\C))^2 \\ 
@V (X_1,X_2)  VV @VV  (X_1,X_2) V\\ 
 (\P^1(\C))^2 @>(A_1,A_2)>> (\P^1(\C))^2 
\end{CD} 
\ee
commutative for some rational function $B$. It follows now from Theorem  \ref{7} that there exist  rational functions 
\be \l{xora1} X_{1,1}, X_{1,2},\dots ,X_{1,l_1} \ \ \ \  {\rm and} \ \ \ \   X_{2,1}, X_{2,2},\dots ,X_{2,l_2}\ee
such that any irreducible invariant curve $\f C$  of bi-degree $(d_2,d_1)$ is parametrized by 
rational functions $X_1$ and $X_2$ satisfying
\be \l{xora2} X_1=X_{1,j_1}\circ \mu_1, \ \ \ \  X_2=X_{2,j_2}\circ \mu_2\ee for some $j_1,$ $1\leq j_1\leq l_1,$ $j_2,$ $1\leq j_2\leq l_2,$ and  $\mu_1, \mu_2 \in \Aut(\P^1(\C))$. Moreover, the numbers   
$l_1$ and $l_2$ can be bounded in terms of 
$d_1,d_2$, and $m$.
Since a parametrization $X_1,X_2$ of $\f C$ is defined in a unique way up to the change 
$$(X_1,X_2)\rightarrow  (X_1\circ \alpha,X_2\circ \alpha), \ \ \ \ \ \alpha\in \Aut(\P^1(\C)),$$ this implies that to prove the theorem it is enough to show that for any fixed rational functions  $X_1,$ $X_2$ 
there exist at most finitely many $\mu\in \Aut(\P^1(\C))$ 
such that the diagram \be \l{opa}
\begin{CD} 
(\P^1(\C))^2 @>(C,C)>>(\P^1(\C))^2 \\ 
@V (X_1 ,X_2\circ \mu )  VV @VV  (X_1 ,X_2\circ \mu) V\\ 
 (\P^1(\C))^2 @>(A_1,A_2)>> (\P^1(\C))^2 
\end{CD} 
\ee
commutes for some rational function $C$, and that the number of such $\mu$ can be bounded in terms of the numbers 
$m,d_1,d_2.$ 

By Lemma \ref{bug}, there exist $B_{1,1},B_{1,2},\dots ,B_{1,s_1}$ and $B_{2,1},B_{2,2},\dots ,B_{2,s_2}$, where $s_1$ and $s_2$ are bounded in terms of 
$m,d_1,d_2$,
such that \eqref{opa} holds if and only if 
$$ C =B_{1,j_1}, \ \ \ \ \ \mu\circ C \circ \mu^{-1}=B_{2,j_2}$$
 for some $j_1,$ $1\leq j_1 \leq s_1$,  $j_2,$ $1\leq j_2 \leq s_2$ and $\mu\in \Aut(\P^1(\C)).$ Thus, we only must show that  for each pair $j_1,$ $j_2$ 
the number of $\mu\in \Aut(\P^1(\C))$ such that 
\be \l{buri} \mu\circ B_{1,j_1} \circ \mu^{-1}=B_{2,j_2}\ee is finite and can be bounded in terms of $m.$ For this purpose, we observe that if 
along with \eqref{buri} the equality $$\t\mu\circ B_{1,j_1} \circ \t\mu^{-1}=B_{2,j_2}$$ holds for some $\t\mu\in \Aut(\P^1(\C))$, then 
 $ \t \mu\circ\mu^{-1}$ belongs to $G_2(B_{2,j_2})$. Therefore, the number of $\mu\in \Aut(\P^1(\C))$ satisfying \eqref{buri} is equal to the order of the group $G_2(B_{2,j_2})$, which is finite by Lemma \ref{kroo}.

Assume finally that both $A_1$ and $A_2$ are Latt\`es maps. In this case, by \linebreak Theorem \ref{in1}  
there exist  a compact Riemann surface   $R$ of genus zero or {\it one}, and 
holomorphic maps $X_1:R\rightarrow \P^1(\C)$ and $X_2:R\rightarrow \P^1(\C)$ of degrees $d_1$ and $d_2$  correspondingly such that  the diagram 
$$
\begin{CD} 
R^2 @>(B,B)>>R^2 \\ 
@V (X_1,X_2)  VV @VV  (X_1,X_2) V\\ 
 (\P^1(\C))^2 @>(A_1,A_2)>> (\P^1(\C))^2 
\end{CD} 
$$
commutes for some holomorphic map $B:R\rightarrow R$. In turn, the commutativity of this diagram implies that for every $d\geq 1$ 
the algebraic curves 
$$A_i^{\circ d}(x)-B(y)=0, \ \ \ \ \ i=1,2,$$ have a factor of genus zero or one. 
By Theorem 3.5 of \cite{aol}, this implies that $X_i$ is a compositional left factor of 
$\theta_{\f O_0^{A_i}}$. Therefore, $\f O_2^{X_i}\preceq \f O_0^{A_i},$ by Lemma \ref{xriak+}.  
Thus,  $\chi(\f O_2^{X_i})\geq 0$, and arguing as in Lemma \ref{lu} we see that, up to the change $$X\rightarrow X\circ \alpha, \ \ \ \alpha\in \Aut(R),$$ there exist only finitely 
many choices for $X_i$. Now we can finish the proof as above using the full versions of 
 Lemma \ref{kroo} and Lemma \ref{bug}. \qed

\end{section}

\vskip 0.2cm

\noindent {\bf Availability of data and material.} The manuscript has no associated data.

\vskip 0.2cm

\noindent {\bf Conflict of interest statement}. The corresponding author states that there is no
conflict of interest.

\end{document}